[1994/12/01]
\documentclass{amsart}
\usepackage{times}
\usepackage{amsmath,amsfonts}
\usepackage{bm}
\usepackage{enumerate}
\usepackage{epsfig}
\numberwithin{equation}{section}

\chardef\bslash=`\\ 

\newtheorem{theorem}{Theorem}[section]
\newtheorem{lemma}[theorem]{Lemma}
\newtheorem{proposition}[theorem]{Proposition}

\theoremstyle{definition}
\newtheorem{definition}[theorem]{Definition}

\theoremstyle{remark}
\newtheorem{remark}[theorem]{Remark}

\numberwithin{equation}{section}

\newcommand\la{\lambda}

\def\un{{\relax\ifmmode 1\!\!1\else$1\!\!1$\fi}}

\newcommand\ho{H\"{o}lder }

\def\ep{\varepsilon}

\def\bs{\backslash}

\def\ra{\rightarrow}

\def\R{\mathbb{R}}

\def\zu{[0,1]}
\def\om{\omega}

\begin{document}




\title[Multifractality and time subordination]{On multifractality and time subordination \\for continuous functions}


\author {St\'ephane Seuret}

\email{seuret@univ-paris12.fr}
    \address{Laboratoire d'Analyse et de Math\'ematiques Appliqu\'ees - Universit\'e Paris-Est - UFR   Sciences et Technologie - 
61, avenue du G\'en\'eral de Gaulle, 94010 Cr\'eteil Cedex, France}


%


%


\begin{abstract} Let $Z:\zu\ra\R$ be a  continuous function.
We show that if $Z$ is "homogeneously multifractal"   (in a sense we precisely define), then $Z$ is the composition of a monofractal function $g$ with a time subordinator $f$ (i.e. $f$ is the integral of a positive Borel measure supported by  $\zu$). When the initial function $Z$ is given, the monofractality exponent of the associated function $g$ is uniquely determined. We study in details a classical example of multifractal functions $Z$, for which we exhibit the associated functions $g$ and $f$. This provides new insights into the understanding of multifractal behaviors of functions.
\end{abstract}

\maketitle

\section{Introduction and motivations}
\label{intro}

Local regularity and multifractal analysis have become unavoidable issues in the past years. Indeed, physical phenomena exhibiting wild local regularity properties have been discovered  in many contexts (turbulence flows, intensity of seismic waves, traffic analysis,..). From a mathematical viewpoint, the multifractal approach is also a fruitful source of interesting problems. Consequently, there is a strong need for a better theoretical understanding of the so-called multifractal behaviors. In this article, we investigate the relations between   multifractal properties and time subordination for continuous functions.

\smallskip

 The most common  functions or processes  used to model irregular phenomena are mono\-fractal, in the sense that they exhibit the same local regularity at each point. Let us recall how the local regularity of a function is measured.
\begin{definition}
Let $Z \in L^\infty_{loc}(\zu)$. For $\alpha\geq 0$ and $t_0\in \zu$, $Z$ is said to belong to $C^\alpha_{t_0}$ if there are a polynomial $P$ of degree less than $[\alpha]$ and a constant $C$ such that, locally around $t_0$,
\begin{equation}
\label{defpoint}
|Z(t) - P(t-t_0)| \leq C |t -t_0|^\alpha.
\end{equation}
The pointwise \ho exponent of $Z$ at $t_0$ is $h_Z(t_0) = \sup\{\alpha\geq 0: \ f\in C^\alpha_{t_0} \}.$

 The singularity spectrum of  $Z$ is then defined by $d_Z(h)=\dim \{t: h_Z(t)=h\}$ ($\dim$ stands for the Hausdorff dimension, and $\dim \emptyset = -\infty$ by convention).
\end{definition}
Hence,  a function $Z:\zu\ra\R$ is said to be monofractal with exponent $H>0$ when  $h_Z(t)=H$ for every $t\in \zu$. For monofractal functions $Z$, $d_Z(H) =1$, while $d_Z(h)=-\infty$ for $h\neq H$. Sample paths of Brownian motions or fractional Brownian motions are known to be almost surely monofractal with exponents less than 1. For reasons that  appear below, {\bf we  focus on monofractal functions associated with an exponent $H \in (0,1]$.}

More complex models had to be used and/or developed, for at least three reasons: the occurrence of intermittence phenomena (mainly in fluid mechanics),   the presence of oscillating patterns (for instance in image processing), or the presence of discontinuities (in finance or telecommunications). Such models may have multifractal properties, in the sense that the support of their singularity spectrum is not reduced to a single point. Among these processes, whose local regularity  varies badly from one point to another,  let us mention Mandelbrot multiplicative cascades and their extensions \cite{BMP,MANDEL2,kahane,BM1} ,  (generalized) multifractional Brownian motions \cite{JLV,BJR} and L\'evy processes \cite{BERTOIN,JAFFLEVY} (for discontinuous phenomena).

Starting from a monofractal process as above in dimension 1, a simple and efficient way to get a more elaborate process is to compose it with a time subordinator, i.e. an increasing function or process. Mandelbrot, for instance, showed the pertinency of   time subordination in the study of financial data \cite{MANDEL}. From a theoretical viewpoint, it is  also challenging to understand how the multifractal properties of a function are modified after a  time change~\cite{RMP,MBLEVY}.

\smallskip

A natural question is to understand the differences between the multifractal processes  above and compositions of monofractal functions with multifractal subordinators.
\begin{definition}
A function  $Z:\zu\ra\R$ is said to be the composition of a monofractal function with a time subordinator (CMT) when $Z$ can be written as
\begin{equation}
\label{decomp}
Z=g\circ f,
\end{equation}
 where $g:\zu\ra\R$ is monofractal with exponent $0<H<1$ and $f:\zu\ra\zu$ is an increasing homeomorphism of $\zu$.
 \end{definition}
In this article, we prove that if   a continuous    function  $Z:\zu\ra\R$ has a "homogeneous multifractal" behavior (in a sense we define just below), then $Z$ is CMT. Hence, $Z$ is the  composition of a monofractal function    with a time subordinator, and   shall  simply be viewed as a complication of a monofractal model. This yields a deeper insight into the understanding of multifractal behaviors of continuous functions, and gives a more important role to the multifractal analysis of  positive Borel measures (which are derivatives of time subordinators).  We explain in Section \ref{self} and \ref{multi} how this decomposition can be used to compute the singularity spectrum of the function $Z$.

%

\medskip

Let us begin with two   cases where a function $Z$ is obviously CMT: 

\smallskip

1. If $Z$ is the  integral of any positive Borel measure $\mu$, then $Z = Id_{\zu}\circ Z$, where the identity $Id_{\zu} $ is   monofractal and $Z$ is increasing. Remark that in this case, $Z$ may even  have exponents greater than 1.

\smallskip

2. Any monofractal function $Z_H$ can be written $Z_H=Z_H \circ Id_{\zu}$, where $ Z_H$   is  monofractal and $ Id_{\zu}$  is undoubtably an homeomorphism of $\zu$. 
 
\smallskip

These two simple cases will be met again below.

\smallskip

To bring   general answers to our problem and thus to exhibit another class of CMT functions , we   develop an approach based on the oscillations of a function $Z:\zu \ra\R$.
For every subinterval $I\subset \zu$, consider the oscillations of order 1 of $Z$ on $I$ defined by
$$\omega_I(Z) =  \sup_{t,t'\in I } |Z(t)-Z(t')|= \sup_{t\in I } Z(t) -\inf_{t\in I } Z(t).$$
 \noindent{\bf In the sequel, we assume that  $Z$ is   continuous and for every non-trivial subinterval $I$ of $\zu$, $\omega_I(Z) >0$. }
 This entails that $Z$  is nowhere locally  constant, which is a natural assumption for the results we are looking for.

 It is very classical that the oscillations of order 1  characterize precisely the pointwise \ho exponents strictly less than 1 (see Section \ref{prel}). 
 

\smallskip

Let us introduce the quantity that will be the basis of our construction.

For every $j\geq 1$, $k\in\{0,..., 2^j-1\}$, we  consider the dyadic intervals $I_{j,k} = [k2^{-j}, (k+1) 2^{-j})$, so that $\bigcup_{k=0,..., 2^j-1} I_{j,k} = [0,1[$, the union being disjoint.
 For every $j\geq 1$ and $k\in\{0,..., 2^j-1\}$, for simplicity we set $\omega_{j,k}(Z) =\omega_{I_{j,k}}(Z) $($=\omega_{\overline{I_{j,k}}}(Z) $ since $Z$ is $C^0$).
 
\begin{definition} 
For every $j\geq 1$, let $H_j(Z)$ be the unique real number such that
\begin{equation}
\label{defhj}
\sum_{k=0}^{2^j-1}( \omega_{j,k}(Z) ) ^{1/H_j(Z)} =1. 
\end{equation}
We then define the  intrinsic monofractal exponent of $Z$ $H(Z)$  as 
\begin{equation}
\label{defh}
H(Z)= \liminf_{j \ra +\infty} H_j(Z).
\end{equation}
\end{definition}
%

This quantity $H(Z)$  characterizes the   asymptotic maximal values of the oscillations    of  $Z$ on the whole interval $\zu$. This exponent is    the core of our theorem, because it gives an upper limit to the maximal time  distortions we are allowed to apply.

It is satisfactory that $H(Z)$ has a functional interpretation.  Indeed, if $Z$ can be decomposed as (\ref{decomp}), then the exponent of the monofractal function $g$ shall not depend on the oscillation approach nor on the dyadic basis. In  Section \ref{secgen} we explain that
\begin{equation}
\label{form}
H(Z)=  \inf\left\{p>0: Z\in B^{1/p,\infty}_{p,{loc}}((0,1))\right\}= \inf\left\{p>0: Z  \in \mathcal{O}_{p}^{1/p}((0,1))\right\} ,
\end{equation}
where $ {B}^{q,\infty}_{1/q,loc}((0,1))$ and $\mathcal{O}_{p}^{1/p}((0,1))$ are respectively  the Besov space and {\em oscillation space} on the open interval $(0,1)$ (see Jaffard in \cite{JAFFBEY} for instance). 

For multifractal functions $Z$ satisfying some multifractal formalism, the exponent $H(Z)$ can also be read on the singularity spectrum of $Z$. Indeed (see Section \ref{secgen}),  $H(Z)$ corresponds to the inverse of the largest possible slope of a straight line going through 0 and tangent to the singularity spectrum $d_Z$ of $Z$.

These remarks are  important to have an idea {\em a priori} of the monofractal exponent of   $g$ in the decomposition $Z=g\circ f$. They also give an intrinsic formula for $H(Z)$.

\medskip

Let us come back to the two simple examples above:  

1. For the integral $Z$ of any positive measure $\mu$,  $\sum_{k=0}^{2^j-1} \omega_{j,k}(Z)   = \sum_{k=0}^{2^j-1} \mu(I_{j,k} ) =1$, hence $H_j(Z)=H(Z)=1$, which corresponds to the monofractal exponent of the identity $Id_{\zu}$ from the oscillations viewpoint. 

\smallskip

2. The first difficulties arise for the monofractal functions $Z_H$. When $Z_H$ is  monofractal of exponent $H$, then we don't have necessarily $H(Z_H)=H$. We always have $H(Z_H)\leq 1$ (see Lemma \ref{lem0} in Section \ref{prel}), but it is always possible to construct wild counter-examples. Nevertheless,  we treat in details the examples of the Weierstrass functions and the sample paths of   (fractional) Brownian motions in Section \ref{mono}, for which the exponent $H(Z_H)$ meets our requirements.

\smallskip 

Unfortunately, the knowledge of $H(Z)$ is not sufficient to get relevant results. 
For instance, consider a function $Z$ that has two different monofractal behaviors on $[0,1/2)$ and $[1/2,1]$. Such an  $Z$ can be obtained as the continuous juxtaposition of two Weierstrass function with distinct exponents $H_1<H_2$: We have $H(Z)=H_1$, and  $Z$ can not be written as the composition of a monofractal function with a time subordinator.  This is a consequence of Lemma \ref{lemmonof},  which asserts that  two monofractal functions $g_1$ and $g_2$ of disctinct exponents $H_1$ and $H_2$ never verify $g_1 = g_2 \circ f$ for any continuous increasing function $f:\zu\ra\zu$ (indeed, such an  $f$ would  "dilate" time everywhere, which is impossible). 


\medskip

We need to introduce a homogeneity condition {\bf C1}  to get rid of these annoying and artificial cases. This condition heuristically imposes that the oscillations of any restriction of $Z$ to a subinterval of $\zu$ have the same asymptotic properties as the oscillations of $Z$ on  $\zu$.
  
\begin{definition}{\bf Condition {\bf C1}:}\\
Let $J\geq 0$, and $K\in\{0,..., 2^j-1\}$. Let $Z_{J,K}$ be   the function 
\begin{eqnarray}
\label{eq00}Z_{J,K} :  
\nonumber t  \in\zu\longmapsto  \frac{Z \circ \varphi_{J,K}(t)}{\omega_{J,K}(Z)} \in \R,
\end{eqnarray}
where $\varphi_{J,K}$ is the canonical affine contraction which maps $\zu$ to $I_{J,K}$.

Condition {\bf C1} is satisfied for $Z$ when there is a real number $H>0$ such that for every $J\geq 0$  and $K\in\{0,..., 2^j-1\}$, $H(Z_{J,K}) =H (= H(Z))$.
\end{definition}
Hence $Z_{J,K}$ is a renormalized version of the restriction of $Z$ to the interval $I_{J,K}$.
Remark that $H(Z_{J,K})$ does not depend on the normalization factor $  {1}/{\omega_{J,K}(Z)}$.
Although self-similar functions are good candidates to satisfy {\bf C1}, a function $Z$ fulfilling this condition does not need at all to possess such a property. 
In order to guarantee that $Z$ is CMT,   we  strengthen the convergence toward $H(Z_{J,K})$.

\begin{definition}
{\bf Condition {\bf C2}:}\\
Assume that Condition {\bf C1} is fulfilled. There are   two positive sequences $(\ep_J)_{J\geq 0}$  and $(\eta_J)_{J\geq 0}$ and two real numbers  $0<\alpha<\beta$ with the following property: 
\begin{enumerate}
\item
$(\ep_J)_{J\geq 0}$  and $(\eta_J)_{J\geq 0}$ are positive non-increasing sequences that converge to zero,
and 
$ \ep_J =o\left(  \frac{1}{ (\log J)^{2+\kappa}}\right)$ for some $\kappa>0$.

%


\medskip

\item
For every $J\geq 0$  and $K\in\{0,..., 2^J-1\}$,  the sequence $(H_j({Z}_{J,K}))_{j\geq 1}$   converges to $H=H({Z }_{J,K})$  (it is not only a liminf, it is a limit) with the following  convergence rate:  
For every $j\geq   [J\eta_J]$,
\begin{eqnarray}
\label{eq1}
 &&|H - H_j({Z }_{J,K}) | \leq \ep_J,\\
 \label{eq1'}
 \mbox{and }  \mbox{for every $k\in\{0,...,2^j-1\}$,}&& 2^ {-j\beta}\leq  \omega_{j,k}(Z_{J,K}) \leq 2^{-j \alpha} .
\end{eqnarray}
\end{enumerate}
\end{definition}
Assuming that $H({Z}_{J,K})$ is a limit is of course a constraint, but not  limiting in practice, since this condition  holds for most of the interesting functions or (almost surely) for most of the sample paths of processes. Similarly, the decreasing behavior 
(\ref{eq1'}) is not very restrictive: such a behavior is somehow expected for a $C^\gamma$ function.

The convergence speed (\ref{eq1}) is a more important constraint, but the convergence rate we impose on $(\ep_J)_{J\geq 0}$ toward 0 is extremely slow, and is realized in the most common cases, as shown below.

\begin{theorem}
\label{maintheo}
Let $Z:\zu \rightarrow \R$ be a  continuous  function.

Assume that $Z$ satisfies {\bf C1} and {\bf C2}.

Then $Z$ is CMT  and the function $g$ in (\ref{decomp}) is  monofractal  of exponent $H(Z)$.
\end{theorem}

\begin{remark}
Such a decomposition is of course not unique: If $Z$ is CMT and $w:\zu\ra\zu$ is $C^\infty$ and strictly increasing, then $Z= (g\circ w) \circ (w^{-1}\circ f)$, where $g\circ w$ is still a monofractal function of exponent $H(Z)<1$ and $w^{-1}\circ f$ is an increasing function.

\smallskip

Nevertheless, if two   decompositions (\ref{decomp}) exist respectively with functions $g_1$, $g_2$, $f_1$ and $f_2$, then $g_1$ and $g_2$ are necessarily monofractal with the same exponent  $H(Z)$. This is again a consequence of Lemma \ref{lemmonof}.
\end{remark}
 
An important consequence of Theorem \ref{maintheo} is that the (possibly) multifractal behavior of $Z$ is contained in the multifractal behavior of $f$. More precisely, since $f$ is an increasing continuous function from $\zu$ to $\zu$, $f$ is the integral of a positive measure, say $\mu$, on $\zu$. The local  regularity of $\mu$ is classically quantified through a local dimension exponent defined for every $t\in\zu$ by
$$\alpha_\mu(t) = \liminf_{r\ra 0^+} \frac{|\log \mu(B(t,r))|}{r} =\liminf_{j\ra+\infty} \frac{|\log_2 \mu(B(t, 2^{-j}))|}{j},$$
where 
$B(t,r)$ stands for the ball (here an interval) with center $t$ and radius $r$, and $|A|$ is the diameter of the set $A$ ($|B(t,r)|=2r$).
The singularity spectrum of $\mu$ is then  
\begin{equation}
\label{defdmu}
\tilde d_\mu(\alpha) = \dim \{t : \alpha_\mu(t)=\alpha\}.
\end{equation}

It is very easy to see that if $\alpha_\mu(t_0)=\alpha$, then $h_f(t_0) = \alpha H$. Hence for every $h\geq 0$, $d_f(h) = \tilde d_\mu(h/H)$, i.e. there is a direct relationship  between the singularity spectrum of $Z$ and the one of $\mu$. As a conclusion, Theorem \ref{maintheo} increases the role of the multifractal analysis of measures, since for the functions satisfying {\bf C1} and {\bf C2}, their multifractal behavior is ruled exclusively by the behavior of $\mu$.
 
\medskip
 
As an application of Theorem \ref{maintheo}, we will prove the following Theorem \ref{thself}, which relates the so-called self-similar functions $Z$ introduced in \cite{JAFFFORM} with the self-similar measures naturally associated with the similitudes defining $Z$.

 Let us recall the definition of self-similar functions. Let $\phi$ be a Lipschitz function on $[0,1]$ (we suppose that the Lipschitz constant $C_\phi$ equals 1, without loss of generality), and let $S_0, S_1, ...., S_{d-1}$ be $d$ contractive similitudes satisfying:
 \begin{enumerate}
 \item
 for every $i\neq j$, $S_i((0,1))\cap S_i((0,1)) =\emptyset $ (open set condition),

\item
 $\displaystyle \bigcup _{i=0}^{d-1} S_i(\zu) = \zu$ (the intervals $S_i(\zu)$ form a covering of $\zu$).

\end{enumerate}
We denote by $0<r_0,r_1,...,r_{d-1}<1$ the ratios of the non trivial similitudes $S_0,...,S_{d-1}$.  By construction $\displaystyle \sum_{k=0 }^{d-1} r_k =1 $.  Let $\la_0, \la_1,...,\la_{d-1}$ be $d$ non-zero real numbers, which satisfy
\begin{equation}
\label{cond1}
0<\chi_{\min} =  \min_{k=0,...,d-1} \left|\frac{r_k}{\la_k}\right| \leq \chi_{\max} =  \max_{k=0,...,d-1} \left|\frac{r_k}{\la_k}\right| <1.
\end{equation}

\begin{definition}
\label{defiself}
A function $Z: \zu\to\zu$ is called self-similar when $Z$ satisfies the following functional equation
\begin{equation}
\label{defself}
\forall \, t\in \zu, \ \ Z(t) = \sum_{k=0}^{d-1} \la_k \cdot (Z\circ (S_k)^{-1}) (t) + \phi(t).
\end{equation}
\end{definition}
Relation (\ref{cond1}) ensures that   $Z$ exists and  is unique~\cite{JAFFFORM}.

Let us consider the unique exponent $\beta>1$ such that 
\begin{equation}
\label{defbeta}
\sum_{k=0} ^{d-1} (\la_k)^\beta =1.
\end{equation} 
This $\beta$ is indeed  greater than 1, since $\sum_{k=0} ^{d-1} r_k =1$ and $|\la_k|>r_k$ for all $k$ by (\ref{cond1}). With the probability vector $(p_0,p_1,..., p_{d-1})=(|\la_0|^\beta , |\la_1|^\beta ,..., |\la_{d-1}|^\beta )$ and the similitudes $(S_k)_{k=0,...,d-1}$ can be associated the unique self-similar   probability measure $\mu$ satisfying
\begin{equation}
\label{defmu}
 \mu = \sum_{k=0} ^{d-1} p_k \cdot  (\mu\circ S_k ^{-1}).
 \end{equation}

\begin{theorem}
\label{thself}
Let $Z$ be defined by (\ref{defself}). Then, either $Z$ is a $\kappa$-Lipschitz function for some constant $\kappa>0$ (expliciteley found in Section \ref{self}), or   $Z$   is CMT and there is a monofractal function $g$ of exponent $1/\beta$ such that
 \begin{equation}
 \label{resthself}
 \mbox{for every $t\in\zu$, } Z(t) =g (\mu[0,t]),
 \end{equation} 
 where $\mu$ is the self-similar measure (\ref{defmu}) naturally associated with the parameters used to define $Z$. 
\end{theorem}

The multifractal analysis of $Z$ follows from the multifractal analysis of $\mu$, which is  a very classical problem (see \cite{BMP}).

\medskip

The paper is organized as follows. In Section \ref{proof}, Theorem \ref{maintheo} is proved, by explicitly constructing the monofractal function $g$ and the time subordinator $f$. Section \ref{secgen} contains the possible extensions of Theorem \ref{maintheo}, the explanation of the heuristics (\ref{form}), and the discussion for exponents greater than 1. In Section \ref{mono}, \ref{self} and \ref{multi}, we detail several classes  of examples to which Theorem \ref{maintheo} applies. First we prove that the usual monofractal functions  $Z$ with exponents $H$ verify {\bf C1} and {\bf C2}. We prove Theorem \ref{thself} in Section \ref{self}. Finally we explicitly compute and plot the time subordinator and the monofractal function     for a classical family of multifractal functions $(Z_a)_{a\in \zu}$ which include Bourbaki's and Perkin's functions.  


 
\medskip
 

Let us finish by the direct by-products and the possible  extensions of this work:

The reader can check that the proof below can be adapted to more general contexts:
\begin{itemize}
\item the dyadic basis can be replaced by any $b$-adic basis.
\item
if  $(\ep_J)$   converges to zero (without any given convergence rate), then (under slight modifications of $(\eta_J)$) the same result holds true. We focused on a simpler case, but in practice, a convergence rate $ \ep_J =o\left(  \frac{1}{ (\log J)^{2+\kappa}}\right)$ shall always be always obtained.
\item
The fact the the quantities $H({Z }_{J,K})$  are limits is only used at the beginning of the proof. In fact, only the existence of  the scale $[J\eta_J]$ such that (\ref{eq1}) and (\ref{eq1'}) hold true at scale $[J\eta_J]$ is determinant. In particular, the conditions may be relaxed: We could treat the case where the $H({Z }_{J,K})$ are only liminf (and not limits). Again,  in practice they are often limits, this is why we adopted this viewpoint.
\end{itemize}

\section{Preliminary results}
\label{prel}

%
%
\subsection{Oscillations and pointwise regularity}

For every $t\in \zu$, let $I_j(t) $ be the unique dyadic interval of generation $j$ that contains $t$, and $I_j^+(t) = I_j(t)+2^{-j}$, $I_j^-(t) = I_j(t)-2^{-j}$. 


Let us recall the characterization of the pointwise \ho exponents smaller than 1 in terms of oscillations of order 1 (see for instance Jaffard in \cite{JAFFBEY}).

\begin{lemma}
\label{lem1}
Let $Z:\zu\ra \R$ a $C^\gamma$ function, for some $\gamma>0$. Assume that $h_Z(t)<1$. Then $$h_Z(t) = \liminf_{r\ra 0^+} \frac{|\log \om_{B(t,r)} (Z)|}{ |\log r|} = \liminf_{j\ra +\infty} \frac{|\log_2 \om_{B(t,2^{-j})}(Z)|}{j}.$$
\end{lemma}

%


In Lemma \ref{lem2}, we impose some uniform behavior of the oscillations of $Z$ on   a nested sequence of coverings of $\zu$. This is used later to prove the monofractality property of the function $g$ in the decomposition $Z=g\circ f$ (Section \ref{proof}),  and also to decompose self-similar functions (in Section \ref{self}).


\begin{lemma}
\label{lem2}
Let $Z:\zu \ra \R$ be a continuous function,  $t\in (0,1)$ and $H\in (0,1)$.

Suppose that there exists an infinite sequence $(T_n)$  of coverings of $\zu$ such that 
\begin{itemize}
\item
each $T_n$ is a finite sequence of disjoint non-trivial intervals of $\zu$, such that  $\bigcup_{T\in T_n} T =\zu$,
\item
$ \lim _{n\ra +\infty} \max_{T\in T_n} |T| =0$,
\item
each interval $T$ in $T_n$ is contained in a unique interval $T'$ of $T_{n-1}$,
\item
for every  $T\in T_n$ and $T\subset T' \in T_{n-1}$, we have $ |T'|^{1+ Z_n } \leq   { |T |} \leq  |T'|  $, for some positive sequence
$(Z_n)$ that converges to zero when $n\ra +\infty$.
\end{itemize}
  Then:
\begin{enumerate}
\item If there exists a positive sequence $(\kappa_n)_{n\geq 1}$ such that for every $T \in T_n$, $\omega_{T}(Z) \leq |T| ^{H-\kappa_n}$, then for every $t\in \zu$, $h_Z(t) \geq H$.
\item
If there exists a positive sequence $(\kappa_n)_{n\geq 1}$ such that for every $T \in T_n$, $\omega_{T}(Z) \geq |T| ^{H+\kappa_n}$, then for every $t\in \zu$, $h_Z(t) \leq H$.
\end{enumerate}
\end{lemma}
Remark that in part (2) of this Lemma, the property needs to be satisfied only  for  a subsequence $({n_k})_{k\geq 1}$ of integers.

\begin{proof}

 Let $t\in (0,1)$, and $r>0$ small enough. For every $n\geq 1$, $t$ belongs to one interval $T\in T_n$, that we denote $T_n(t)$. 
 Denote by $n_r$ the smallest integer $n$ so that $T_n(t) \subset B(t,r)$.  By construction, $t\in T_{n_r-1}(t)$ and $|T_{n_r-1}(t)| \geq r$ (since $ T_{n_r-1}(t)\not\subset B(t,r)$). By the fourth property of the sequence $(T_n)$, we have $ 2r \geq | T_{n_r}(t)| \geq  |T_{n_r-1}(t)|   ^{1+Z_{n_r}} \geq r^{1+Z_{n_r}}$.
  
 Let us start by part (2), which is very easy to get.
We have $\om_{B(t,r)}(Z) \geq \om_{T_{n_r}(t)}(Z) \geq |T_{n_r}(t)|^{H+\kappa_{n_r}} \geq r^{(1+Z_{n_r})(H+\kappa_{n_r})} $.

Applying Lemma \ref{lem1}, and using that $Z_{n_r}$ and $\kappa_{n_r}$ go to zero when $r$ goes to zero, we obtain $h_Z(t) \leq H$.

We now focus on part (1), which is slightly more delicate. If $B(t,r) \subset T_{n_r-1}(t)$, then we have $\om_{B(t,r)} (Z)\leq \om _{T_{n_r-1}(t)}(Z) \leq |T_{n_r-1}(t)|^{ H-\kappa_{n_r-1}} \leq (2r)^{(H-\kappa_{n_r-1})/(1+Z_{n_r-1})}$.
 
 If $B(t,r) \not\subset T_{n_r-1}(t)$, then there is an integer $p$ (which depends on $r$) such that $ B(t,r)\bs T_{n_r-1}(t) $ is covered by one interval $T\in T_{p}$ and not covered by any interval of $T_{p+1}$. Using the same arguments as above, we get  $|T| \leq | B(t,r)\bs T_{n_r-1}(t) | ^{1/(1+Z_{p+1})} \leq r^{1/(1+Z_{p+1})}$ (remark that $|B(t,r)\bs T_{n_r-1}(t)|\leq r$).

 Now we have
 \begin{eqnarray*}
 \om_{B(t,r)} (Z)& \leq &\om _{T_{n_r-1}(t)}(Z)+\om _{B(t,r)\bs T_{n_r-1}(t)}(Z) \\
 & \leq &(2r)^{(H-\kappa_{n_r-1})/(1+Z_{n_r-1})} + |T|^{H-\kappa_{p}}\\
 & \leq &(2r)^{(H-\kappa_{n_r-1})/(1+Z_{n_r-1})} + r^{(H-\kappa_p)/(1+Z_{p+1})}
 \end{eqnarray*}
 Since $\kappa_{n_r}$, $Z_{n_r}$, $\kappa_p$ and $Z_{p}$ converge to 0 as $r \ra 0$, Lemma \ref{lem1} yields $h_Z(t) \geq H$.
\end{proof}
%



\subsection {Two easy properties for the study of $H(Z)$}

Let us begin with an easy upper-bound for $H(Z)$.
\begin{lemma}
\label{lem0}
Let $Z:\zu \ra \R$ be a non-constant continuous function. Then $H(Z)\leq 1$.
\end{lemma}
 \begin{proof}
 We can assume without loss of generality that $\om_{\zu}(Z) =1$.
 Let $j\geq 1$.
 By construction, $\sum_{k=0}^{2^j-1} \om_{j,k}(Z) \geq 1$. In order to have (\ref{defhj}), we necessarily have $H_j(Z)\leq 1$. Hence the result.
\end{proof}

\begin{lemma}
\label{lemmonof}
Let $g_1$ and $g_2$ be two real monofractal functions on $\zu$ of disctinct exponents $0<H_1 < H_2<1$. There is no continuous strictly increasing function $f:\zu\ra\zu$ such that $g_1 = g_2 \circ f$.
\end{lemma}
 \begin{proof}
 Suppose that such a function $f$ exists. Let $\ep>0$. This function $f$ is  Lebesgue-almost everywhere differentiable. There is  a set $E$ of positive Lebesgue measure such that for every $t\in E$, $f'(t) >0$.
 Around such a $t$, we have $f(t+h) - f(t) = f'(t)h +o(t)$. Consequently, since $h_{g_2}(f(t)) = H_2$, for every $|h|$ small enough we have
 $$ |(g_2\circ f)(t+h) - (g_2\circ f) (t)| \leq |f(t+h) - f(t) | ^{H_2 -\ep} \leq C |h|^{H_2 -\ep} .$$
 This shows that $h_{g_2\circ f} (t)\geq H_2$.
 Using again that $h_{g_2}(f(t)) = H_2$,  there is a sequence $(h'_n)_{n\geq 1}$ converging to zero such that for every $n\geq 1$, $|g_2(f(t)+ h'_n) - g_2(f(t))| \geq  |h'_n|^{H_2+\ep}$. Choosing $h_n$ so that $f(t+h_n) = f(t)+h'_n$, we see that
 $$| (g_2\circ f)(t+h_n) - (g_2\circ f) (t)| \geq |f(t+h_n) - f(t) | ^{H_2 +\ep} \geq C |h_n|^{H_2 +\ep} .$$
This holds for an infinite number of real numbers $(h_n)$ converging to zero. Hence $h_{g_2\circ f} (t)= H_2$, which contradicts $h_{g_1}(t) = H_1$. 
 \end{proof}

\subsection {A functional interpretation of $H(Z)$}

Note first that the previous results hold in the case where a $b$-adic basis, $b\geq 2$, is used instead of the dyadic basis.
In fact, there is a functional interpretation of the exponent $H(Z)$, independent of any basis, provided by the {\em Oscillation spaces} of Jaffard \cite{JAFFBEY} and the Besov spaces. Let us recall their definition, that we adapt to our context of nowhere differentiable functions.

Let $Z$ be a $C^\gamma$ function on $(0,1)$, where $C^\gamma$ is the global homogeneous \ho space and $\gamma>0$. Since \cite{JAFFFORM} where the theoretical foundations of multifractal analysis of functions were given, a  quantity classically considered when performing the  multifractal analysis of $Z$ is the scaling function $\eta_Z(p)= \sup\left \{s>0: Z\in B^{s/p,\infty}_{p,{loc}}((0,1))\right\}$. 

Later, in \cite{JAFFBEY}, Jaffard also proved the pertinency   in multifractal analysis of his oscillation spaces $\mathcal{O}^{s/p}_{p}((0,1))$, whose definitions are based on wavelet leaders (we do not need much more details here). He also considered the associated scaling function $\zeta_Z(p) = \sup\left\{s>0:  Z \in \mathcal{O}^{s/p}_{p}((0,1))\right\}$.

Finally, still in \cite{JAFFBEY}, Jaffard   studied the spaces $\mathcal{V}^{s/p}_p((0,1))$, which are closely related to our exponent $H(Z)$, defined as follows: Denote, for $j\geq 1$ and $k\in \{0,...,2^j-1\}$, $\Omega_{j,k} (Z)= \omega_{[ k2^{-j} - 3 2^{-j}, k2^{-j} + 3 2^{-j} ]}(Z)$, and consider the associated scaling function (we assume hereafter that $Z$ is nowhere differentiable, as in Theorem \ref{maintheo})
$$\nu _Z (p) = 1+ \liminf_{j\ra +\infty} \frac{\log_2 \sum_{k=0}^{2^j-1} (\Omega_{j,k} (Z))^p } { -j}.$$
For $p>0$ fixed, it is obvious that there is a constant $C_p>1$ such that 
$$1/C_p\sum_{k=0}^{2^j-1} (\omega_{j,k} (Z))^p \leq \sum_{k=0}^{2^j-1} (\Omega_{j,k} (Z))^p \leq C_p \sum_{k=0}^{2^j-1} (\omega_{j,k} (Z))^p,$$
since $ \omega_{j,k}(Z) \leq \Omega_{j,k}(Z) \leq \sum_{l\in\{-3,-2,...,2,3\}} \omega_{j,k+l}(Z) $. As a consequence, $\nu _Z(p) = 1+ \liminf_{j\ra +\infty} \frac{\log_2 \sum_{k=0}^{2^j-1} (\omega_{j,k} (Z))^p } { -j}.$ Comparing   the definition of $H(Z)$ with this formula, we easily see that $H(Z)$ is the unique positive real number such that $\nu_Z(1/H(Z)) =1$. 

\smallskip

The main point is that the three scaling functions $\eta_Z$, $\zeta_Z$ and $\nu _Z$ coincide as soon as $p\geq 1$ \cite{JAFFBEY}, and $\eta_Z (1/H(Z)) = \zeta_Z (1/H(Z))=1$. Using the property of the Besov domains,  we have
$$H(Z) ^{-1} =  \inf\left\{p>0: Z  \in B^{1/p,\infty}_{p,\mbox{loc}} ((0,1))\right\}= \inf\left\{p>0: Z   \in \mathcal{O}^{1/p}_{p} ((0,1))\right\} .$$

\subsection{Precisions for functions satisfying a multifractal formalism}

Consider the scaling function $\zeta_Z(p)$ above. 
Then for any function $Z$ having some global \ho regularity \cite{JAFFBEY}, $Z$ is said to obey the multifractal formalism for functions if its singularity spectrum is obtained as the Legendre transform of its scaling function, i.e.
$$\mbox{for every $h\geq 0$, } \ \ d_Z(h) = \inf _{p \in \R} (ph -\zeta_Z(p)+1) \ \ (\in \R^+ \cup\{-\infty\}) .
$$
In particular, since $\zeta_Z(1/H(Z))=1$, we always have $d_Z(h) \leq h/ H(Z)$  (by using $p=1/H(Z)$ in the inequality above). 

Moreover,   assume that $h_c=\zeta'_Z(1/H(Z))$ exists and that $Z$ satisfies the multifractal formalism associated with $\zeta_Z$ at the exponent $h_c$. This means that the  inequality above holds true for   $h=h_c$, i.e.  $d_Z(h_c) = h_c/H(Z)$. 

From the two last properties we get that $1/H(Z)$ is the slope of the tangent to the (concave hull of the) singularity spectrum of $Z$, as claimed in the introduction.

\section{Proof of the decomposition of Theorem \ref{maintheo}}
\label{proof}

The functions $g$ and $F$ are constructed iteratively. First remark that since $(\eta_j)$ converges to zero, one can also assume, by first replacing $\eta_j$ by $\max(\eta_j, 1/\log j)$ and then by imposing that $(\eta_j)$ is non-increasing, that the   sequence  $(\eta_j)$ satisfies:
\begin{itemize}
\item for every $j\geq 1$, $\eta_j \geq 1/\log j$,
\item
$(j\eta_j)$ is now a non-decreasing sequence and $j\eta_j \ra +\infty$ when $j\ra +\infty$,
\item  $(\eta_j)$ still satisfies  (\ref{eq1}) and (\ref{eq1'}).
\end{itemize}
Assume that conditions {\bf C1} and {\bf C2} are fulfilled.

\subsection{First step  of the construction of  $g$ and $f$}

The exponent $H(Z_{0,0})=H(Z)=H$ is the limit of the sequence $H_j(Z)$, so there exists  a generation $J_0\geq 1$ such that for every $j\geq J_0$, $ |H - H_{ j}(Z) | \leq \ep_{0}.$

We set  $H_0=H_{J_0}(Z)$, and by construction  we  have 
 $\sum_{k=0}^{2^{J_0}-1}( \omega_{J_0,k}(Z) ) ^{1/H_0} =1.$

\medskip

We then define the first step of the construction of the function $f$: we set \begin{eqnarray*}
f_0(t) = \sum_{ k'=0}^{k-1} ( \omega_{J_0,k'}(Z) ) ^{1/H_0} +  ( \omega_{J_0,k }(Z) ) ^{1/H_0} (2^{J_0}
t-k) \ \mbox{ if } t\in I_{J_0,k}.
\end{eqnarray*}
 This function $f_0$ is strictly increasing, continuous and  affine on each dyadic interval. Moreover,     $f_0(\zu) = \zu$. Let us denote $U_{J_0,k}$ the image of the interval $I_{J_0,k}$ by $f_0$, for every $k\in\{0,...,2^{J_0}-1\}$. The set of intervals $\{U_{J_0,k}: k \in\{0,...,2^{J_0}-1\}$ clearly forms a partition of $[0,1)$. One remarks that 
\begin{equation}
\label{eq3}
\forall k\in\{0,...,2^{J_0}-1\},  \ \ |f_0 (I_{J_0,k} ) | = |U_{J_0,k}| = ( \omega_{J_0,k }(Z) ) ^{1/H_0}  .
\end{equation}

\medskip

The first step of the construction of $g$ is then naturally achieved as follows: we set
\begin{eqnarray*}
g_0(y) & = & Z((f_0)^{-1}(y))  \ \mbox {for $y\in \zu$},\\
\mbox {or equivalently } \ g_0(f_0(t)) & = &  Z(t)  \ \mbox {for $t\in \zu$}.
\end{eqnarray*}
 This function $g_0$ maps any interval $U_{J_0,k}$ to the interval $Z(I_{J_0,k})$, and thus satisfies:
 $$ \omega _{U_{J_0,k} }({g_0}) = \omega_{J_0,k } (Z)=  |U_{J_0,k}| ^{H_0}.$$

As a last remark, there are two real numbers $0<\alpha'<\beta'$ such that for every $k$ $2^{-J_0\beta'/H_0} \leq  |U_{J_0,k}| \leq 2^{-J_0\alpha'/H_0}$. Without loss of generality, we can assume that $\alpha' =\alpha$ and $\beta'=\beta$ ($\alpha$ and $\beta$ appear in condition {\bf C2}) by changing $\alpha$ into $\min(\alpha',\alpha)$ and $\beta=\max(\beta,\beta')$, so that  
\begin{equation}
\label{eq22}
\mbox{for every $k$, } \ 2^{-J_0\beta/H_0} \leq  |U_{J_0,k}| \leq 2^{-J_0\alpha/H_0}.
\end{equation}

\subsection{First iteration to get the second step of  the construction of  $g$ and $f$}

We perform the second step of the construction. Let us focus on one interval $I_{J_0,K}$, on which we refine the behavior of $f_0$. By condition {\bf C2} and especially  (\ref{eq1}), we have
\begin{equation}
\label{eq5}
\sum_{k'=0,...,2^{ [J_0\eta_{J_0}]}-1} (\omega_{ [J_0\eta_{J_0}],'k} (Z_{J_0,K}))^{1/H_1} =1,
\end{equation}
where $H_1 = H_{[J_0\eta_{J_0}]}(Z_{J_0,K})$ satisfies $\ |H - H_1| \leq \ep_{J_0}$. 

Let $J_1= J_0+ [J_0\eta_{J_0}]$, hence $J_1-J_0= [J_0\eta_{J_0}]$.  Remark that, by (\ref{eq1'}), we have  for every $k'\in\{0,...,2^{ [J_0\eta_{J_0}]}-1\}$
\begin{equation}
\label{eq20}
 2^{- [J_0\eta_{J_0}] \beta} \leq | \omega_{ [J_0\eta_{J_0}],k'} (Z_{J_0,K}) | \leq  2^{- [J_0\eta_{J_0}] \alpha}.
 \end{equation}
\noindent Now, remembering the definition of $Z_{J_0,K}$, we obtain that for every $k'\in\{ 0,...,2^{J_1-J_0}-1\}$,
\begin{equation}
\label{eq21}
\om_{ J_1-J_0,k'} (Z_{J_0,K}) = \frac{  \om_{J_1,K 2^{J_1-J_0}+k'}(Z)}{\om_{J_0,K}(Z)}.
 \end{equation}
  Consequently, (\ref{eq5})  is equivalent to
$$\sum_{k=0,...,2^{J_1}-1: I_{J_1,k}\subset I_{J_0,K}} (\omega_{ J_1,k} (Z))^{1/H_1} =(\omega_{ J_0,K} (Z))^{1/H_1} ,$$
 and thus
 $$\sum_{k=0,...,2^{J_1}-1: I_{J_1,k}\subset I_{J_0,K}} (\omega_{ J_1,k} (Z))^{1/H_1} (\omega_{ J_0,K} (Z))^{1/H_0-1/H_1}  =(\omega_{ J_0,K} (Z))^{1/H_0} .$$

We now define the function $f_1$ as a refinement on $f_0$ on the dyadic interval $I_{J_0,K}$. We set  for every $ k\in\{K2^{J_1-J_0},..., (K+1)2^{J_1-J_0}-1\}$ and for $t\in I_{J_1,k}$
\begin{eqnarray*}
f_1(t) & = & f_0(K2^{-J_0}) \\ &+&\sum_{ k'=K2^{J_1-J_0}}^{k-1} ( \omega_{J_1,k'}(Z) ) ^{1/H_1} (\omega_{ J_0,K} (Z))^{1/H_0-1/H_1} \\ &+&  ( \omega_{J_1,k}(Z) ) ^{1/H_1}(\omega_{ J_0,K} (Z))^{1/H_0-1/H_1}  (2^{J_1} t-k) .
\end{eqnarray*}

This can be achieved simultaneously  on every dyadic interval $I_{J_0,K}$, $K\in\{0,...,2^{J_0}-1\}$, by using the same generation $J_1$ for the subdivision (indeed, condition {\bf C2} ensures that the convergence rate of $H_j(Z_{J_0,k})$ does not depend on $k$). The obtained function is again an increasing continuous function,  affine on every dyadic interval of generation $J_1$.

Let us denote $U_{J_1,k}$ the image of the interval $I_{J_1,k}$ by $f_1$, for every $k\in\{0,...,2^{J_1}-1\}$. The set of intervals $\{U_{J_1,k}: k \in\{0,...,2^{J_1}-1\}$ again forms a partition of $[0,1)$. We get
\begin{equation}
\label{eq4}
\forall k\in\{0,...,2^{J_1}-1\},  \ \ |U_{J_1,k}| = ( \omega_{J_1,k }(Z) ) ^{1/H_1}(\omega_{ J_0,K} (Z))^{1/H_0-1/H_1} .
\end{equation}
but the main point is that we did not change the size of the oscillations of $f_0$ on the dyadic intervals of generation $J_0$, i.e.  $f_1(I_{J_0,K}) = f_0(I_{J_0,K}) $.

\medskip

The second step of the construction of $g$ is realized by refining the behavior of $g_0$: Set
\begin{eqnarray*}
g_1(y) & = & Z((f_1)^{-1}(y))  \ \mbox {for $y\in \zu$}.
\end{eqnarray*}
 This function $g_1$ maps any interval $U_{J_1,k}$ to the interval $Z(I_{J_1,k})$, and thus satisfies:
 $$ \omega _{U_{J_1,k} }({g_1}) = \omega_{J_1,k } (Z) \mbox { with } |U_{J_1,k}| = ( \omega_{J_1,k }(Z) ) ^{1/H_1}(\omega_{ J_0,K} (Z))^{1/H_0-1/H_1} .$$

Finally, we want to compare the size  of the interval $U_{J_1,k}$ with the size of its father interval (in the preceding generation) $U_{J_0,K}$. For this, let us choose $k\in\{ 0,...,2^{J_1 }-1\}$ and $K\in\{ 0,...,2^{ J_0}-1\}$ are such that  $I_{J_1,k}\subset I_{J_0,K}$ (hence $k$ can be written $k=K.2^{J_1-J_0} +k'$ with $k'\in\{0,...,2^{J_1-J_0}-1\}$). Then, by  (\ref{eq21}),
\begin{eqnarray*}
 |U_{J_1,k}| & = &  ( \omega_{J_0,K }(Z))^{1/H_0}  (\omega_{J_1-J_0,k' }(Z_{J_0,K}) ) ^{1/H_1}\\
  & = &  |U_{J_0,k}|(\omega_{J_1-J_0,k' }(Z_{J_0,K}) ) ^{1/H_1}.
 \end{eqnarray*}
Using  (\ref{eq20}) we get
\begin{eqnarray*}
 |U_{J_1,k}|    \geq      |U_{J_0,k}| 2^{-(J_1-J_0)\beta/H_1}= |U_{J_0,k}| 2^{-[J_0\eta_{J_0}]\beta/H_1}.
 \end{eqnarray*}
On the other side, we know by (\ref{eq22}) that $ |U_{J_0,k}| \leq 2^{-J_0 \alpha/H_0}$, hence
\begin{eqnarray*}
 |U_{J_0,k}| \geq  |U_{J_1,k}|    \geq      |U_{J_0,k}|^{1+ \eta_{J_0}\frac{\beta H_0}{\alpha H_1}} ,
 \end{eqnarray*}
where the left inequality simply comes from the fact that $I_{J_1,k}\subset I_{J_0,K}$.

\subsection{General iterating construction of  $g$ and $f$}

This procedure can be iterated. Assume that the sequences $(J_p)_{p\geq 1}$, $(f_p)_{p\geq 1}$ and $(g_p)_{p\geq 1}$ are constructed for every $p\leq n$, and that they satisfy:
\smallskip

\begin{enumerate}
\item 
for every $1\leq p\leq n$, $J_p=J_{p-1} +[J_{p-1}\eta_{J_{p-1}}]$ and $|H-H_p| \leq \ep_{J_{p-1}}$,
\smallskip

\item
for every $1\leq p\leq n$, $f_p$ is a continuous strictly increasing function, affine on each dyadic interval $I_{J_p,k}$ and if we set $f_p(I_{J_p,k}) = U_{J_p,k}$, then
\begin{equation}
\label{eq7}
|f_p(I_{J_p,k})| = |U_{J_p,k}| = ( \omega_{J_p,k }(Z) ) ^{1/H_p} \prod_{m=0}^{p-1}(\omega_{ J_m,K_m(k)} (Z))^{1/H_m-1/H_{m+1}},
\end{equation}
where $K_m(k)$ is the unique integer such that $I_{J_p,k } \subset I_{ J_m,K_m(k)}$, for $m<p$,
\smallskip

\item
For every $1 \leq p\leq n$, the set of intervals $\{U_{J_p,k}: k \in\{0,...,2^{J_p}-1\}$ forms a partition of $[0,1)$. 
\smallskip

\item
For every $1\leq p\leq n$, if $U_{J_n,k} \subset U_{J_{p-1},K_{p-1}(k)}$, then $$ {|U_{J_{p-1},K_{p-1}(k)}|} ^{1+ \eta_{J_{p-1}}\frac{\beta H_{p-1}}{\alpha H_{p}}}  \leq {|U_{J_p,k}|}\leq {|U_{J_{p-1},K_{p-1}(k)}|} ,$$
\smallskip

\item
for every $1\leq p\leq n$,  for $y\in \zu$, $g_p(y) = Z((f_p)^{-1}(y)) $
\smallskip

\item
for every $1\leq p\leq n$,   for every $k\in\{0,...,2^p-1\}$, we have $f_m(k2^{-p}) = f_p(k2^{-p}) $ for every $p\leq m\leq n$.
\end{enumerate}
\smallskip

The last item ensures that once the value of $ f_p$ at $k2^{-p}$ has been chosen, every $f_m$, $m\geq p$, will take the same value  at $k2^{-p}$.

\medskip

To build $f_{n+1}$ and $g_{n+1}$, the procedure is the same as above. We use $J_{n+1} = J_n + [J_n\eta_{J_n}]$, and we focus on one interval $I_{J_n,K}$.
We have by   (\ref{eq1}) \begin{eqnarray*} 
\sum_{k=0,...,2^{J_{n+1}-J_n}-1 } (\omega_{ J_{n+1}-J_n,k} (Z_{J_n,K}))^{1/H_{n+1} }   =( \omega_{J_n,K }(Z) ) ^{1/H_{n+1}}  ,
\end{eqnarray*}
where $H_{n+1} = H_{[J_n\eta_{J_n}]}(Z_{J_n,K})$ satisfies $\ |H - H_{n+1}| \leq \ep_{J_n}$. we have  for every $k'\in\{0,...,2^{ [J_n\eta_{J_n}]}-1\}$
\begin{equation}
\label{eq30}
 2^{- [J_n\eta_{J_n}] \beta} \leq | \omega_{ [J_n\eta_{J_n}],k'} (Z_{J_n,K}) | \leq  2^{- [J_n\eta_{J_n}] \alpha}.
 \end{equation}
and
\begin{equation}
\label{eq31}
\om_{ J_{n+1}-J_n,k'} (Z_{J_n,K}) = \frac{  \om_{J_{n+1},K 2^{J_{n+1}-J_n}+k'}(Z)}{\om_{J_n,K}(Z)}.
 \end{equation}
The same manipulations as above yield
\begin{eqnarray} 
\label{eq8}
&&\sum_{k=0,...,2^{J_{n+1}}-1: I_{J_{n+1},k}\subset I_{J_n,K} } (\omega_{ J_{n+1},k} (Z))^{1/H_{n+1} }\prod_{m=0}^{n}(\omega_{ J_m,K_m(k)} (Z))^{1/H_m-1/H_{m+1}} \\ \nonumber&& =( \omega_{J_n,K }(Z) ) ^{1/H_n} \prod_{m=0}^{n-1}(\omega_{ J_m,K_m(k)} (Z))^{1/H_m-1/H_{m+1}} ,
\end{eqnarray}

Then $f_{n+1}$ is a refinement on $f_n$: For every $ k\in\{K2^{J_{n+1}-J_n},..., (K+1)2^{J_{n+1}-J_n}-1\}$ and for $t\in I_{J_{n+1},k}$
\begin{eqnarray*}
f_{n+1}(t) &=  &  f_n(K2^{-J_n}) \\ & +&\sum_{ k'=K2^{J_{n+1}-J_n}}^{k-1} ( \omega_{J_{n+1},k'}(Z) ) ^{1/H_{n+1}} \prod_{m=0}^{n}(\omega_{ J_m,K_m(k')} (Z))^{1/H_m-1/H_{m+1}}  \\ & + &\  ( \omega_{J_{n+1},k}(Z) ) ^{1/H_{n+1}}\prod_{m=0}^{n}(\omega_{ J_m,K_m(k)} (Z))^{1/H_m-1/H_{m+1}}  (2^{J_{n+1}} t-k) .
\end{eqnarray*}
Remark that for every $(k,k') \in\{K2^{J_{n+1}-J_n},..., (K+1)2^{J_{n+1}-J_n}-1\}^2$, for every $m\in \{0,...., n\}$, $K_m(k)=K_m(k')$.

This can be achieved simultaneously  on every dyadic interval $I_{J_n,K}$, $K\in\{0,...,2^{J_n}-1\}$, by using the same generation $J_{n+1}$ for the subdivision. The obtained function is again an increasing continuous function which is affine on every dyadic interval of generation $J_{n+1}$.  

We then define $g_{n+1}$ by $
g_{n+1}(y) = Z((f_{n+1})^{-1}(y))$ for $y\in \zu$.  Let   $U_{J_{n+1},k}$ be the image of the interval $I_{J_{n+1},k}$ by $f_{n+1}$, for every $k\in\{0,...,2^{J_{n+1}}-1\}$.
 This function $g_{n+1}$ maps any interval $U_{J_{n+1},k}$ to the interval $Z(I_{J_{n+1},k})$, and thus satisfies:
 $$ \omega _{U_{J_{n+1},k} }({g_{n+1}}) = \omega_{J_{n+1},k } (Z)$$
 with
 $$
  |U _{J_{n+1},k}| =( \omega_{J_{n+1},k}(Z) ) ^{1/H_{n+1}}\prod_{m=0}^{n}(\omega_{ J_m,K_m(k)} (Z))^{1/H_m-1/H_{m+1}}  .$$ 
  
  \smallskip
  
At this point, all the items of the iteration are ensured, except the item (4).  We prove it now. As above, let us choose $k\in\{ 0,...,2^{J_ {n+1} }-1\}$ and $K\in\{ 0,...,2^{ J_2}-1\}$ are such that  $I_{J_{n+1},k}\subset I_{J_n,K}$, and let  $k'\in\{0,...,2^{J_{n+1}-J_n}-1\}$ be such that  $k=K.2^{J_{n+1}-J_n} +k'$. We have  by  (\ref{eq31}) 
\begin{eqnarray}
 \nonumber |U_{J_{n+1},k}| & = &  ( \omega_{J_n,K }(Z))^{1/H_0}  (\omega_{J_{n+1}-J_n,k' }(Z_{J_n,K}) ) ^{1/H_{n+1}}\\
 \label{eq40}& = &  |U_{J_{n},k}| (\omega_{J_{n+1}-J_n,k' }(Z_{J_n,K}) ) ^{1/H_{n+1}}.
 \end{eqnarray}
Then, by (\ref{eq20}),  
\begin{eqnarray*}
 |U_{J_{n+1},k}|    \geq      |U_{J_n,k}| 2^{-(J_{n+1}-J_n)\beta/H_n}= |U_{J_n,k}| 2^{-[J_n\eta_{J_n}]\beta/H_{n+1}}.
 \end{eqnarray*}
As above, since by (\ref{eq22}) that $ |U_{J_n,k}| \leq 2^{-J_n \alpha/H_n}$, we have
\begin{eqnarray*}
 |U_{J_n,k}| \geq  |U_{J_{n+1},k}|    \geq      |U_{J_n,k}|^{1+ \eta_{J_n}\frac{\beta H_n}{\alpha H_{n+1}}}.
 \end{eqnarray*}

\subsection{Convergence   of  $(g_n)_{n\geq 0}$ and $(f_n)_{n\geq 0}$}
  \label{secconv}
  
The convergence of the sequence $(f_n)$ to a  function $f$ is almost immediate. Indeed, each $f_n$ is an increasing function from $\zu$ to $\zu$, and by item (5) of the iteration procedure, for every $j\geq 1$, for every $k\in\{0,...,2^j-1\}$, $f_m(k 2^{-j})$ is constant as soon as $J_m \geq j$.  

Recall that for every $m$ and $k$, $|f_m(I_{J_m,k})| = |U_{J_m,k}|$. By (\ref{eq40}), and using (\ref{eq20}), we obtain that $|U_{J_{m+1},k}|\leq |U_{J_m,k}| 2^{-(J_{m+1}-J_m)\alpha},$ and iteratively 
\begin{equation}
\label{tn}
\mbox{for every $m\geq 1$, }\ |U_{J_{m},k}|\leq   C2^{-J_{m}\alpha},
\end{equation}
 for some constant $C$. Hence the sequence $(|U_{J_{m+1},k}|)_{m\geq 1}$ converge exponentially fast  to zero, with an upper bound independent of $k$.

As a consequence, if $m\geq n $, then 
\begin{eqnarray*}
\|f_n-f_m\|_\infty & \leq & \max_{k\in\{0,..., 2^{J_m}-1\} } |f_m(I_{{J_m,k}})| \leq  \max_{k\in\{0,..., 2^{J_m}-1\} }  |U_{J_m,k}|\\
& \leq & C2^{-J_{m }\alpha}.
\end{eqnarray*}

This Cauchy criterion  immediately gives the uniform convergence of the function series $(f_n)$ to a continuous function $f$, whose value at each dyadic number  is known as explained just above. The limit function $f$ is also strictly increasing, since it is strictly increasing on the dyadic numbers.

\medskip

The convergence of the functions sequence $(g_n)_{n\geq 0}$ is then straightforward. Indeed, each $f_n$  is an homeomorphism of $\zu$, and admits a continuous inverse $f_n^{-1}$. We thus have, for every $n\geq 1$, $g_n = Z \circ f_n^{-1}$. The series $(f_n^{-1})$ also converges uniformly on $\zu$. Since $Z$ is uniformly continuous on $\zu$, $(g_n)$ converges uniformly to a continuous function $g:\zu\ra\zu$. 

\medskip

Remark that $f$ also admits an inverse function $f^{-1}$, and that $g=Z\circ f^{-1}$.

\subsection{Properties of  $g$ and $f$}
  
Obviously, $f$ is a strictly increasing function from $\zu$ to $\zu$, which is what we were looking for. All we have to prove the monofractality property of $g$. This will follow from Lemma \ref{lem2}.

\medskip

It has been noticed before that if we set, for every $n\geq 0$, $T_n = \{U_{J_n,k}: k\in \{0,...,2^{J_n}-1\} \}$, then every $T_n$ forms a covering of $\zu$ constituted by  pairwise  distinct intervals. We obviously have:
\begin{itemize}
\item
 $\lim_{n\ra +\infty} \max_{T\in T_n} |T| =0$ (using the remarks of Section \ref{secconv} above),
 \item $(T_n)$ is a nested sequence of intervals,
 \item
 by item (4) of the iteration procedure, if $T \in T_n \subset T'\in T_{n-1}$, then we have ${|T'|} ^{1+ Z_n }  \leq {|T|} \leq {|T'|}$, , with $Z_{n} = \eta_{J_{n-1}}\frac{\beta H_{n-1}}{  \alpha H_{n}}$. This sequence $(Z_n)$ converges to zero, since $(\eta_n)$ converges to zero and $(H_n)$ converges to $H$.
\end{itemize}

In order to apply Lemma \ref{lem2} and to get the monofractality property of $g$, it is thus enough to prove the last required  properties, i.e.   there is a positive sequence $(\kappa_n)$ converging to zero such that for every $T\in T_n$, $ |T|^{H+\kappa_n} \leq \om_{T} (g) \leq |T|^{H- \kappa_n}$.

\medskip

For this, let $n\geq 1$ and $T\in T_n$. This interval $T$ can be written $U_{J_n,k}$ for some $k\in \{0,...,2^{J_n}-1\}$. We have  $|U_{J_n,k}| =  ( \omega_{J_n,k }(Z) ) ^{1/H_n} \prod_{m=0}^{n-1}(\omega_{ J_m,K_m(k)} (Z))^{1/H_m-1/H_{m+1}}$ by construction , and $g(U_{J_n,k})=g_n(U_{J_n,k})=\om_{J_n,k}(Z)$. We just have to verify that  $ |U_{J_n,k}|^{H+\kappa_n}\leq \om_{J_n,k}(Z) \leq |U_{J_n,k}|^{H-\kappa_n}$, for some $\kappa_n>0$ independent of $k$.

We have
\begin{eqnarray*}
\log |U_{J_n,k}| = \frac{1}{H_n} \log \omega_{J_n,k }(Z) + \sum_{m=0}^{n-1}( \frac{1}{H_m} -\frac{1}{H_{m+1}}) \log \omega_{ J_m,K_m(k)} (Z) 
\end{eqnarray*}
Writing that $ \left|\frac{1}{H_m} -\frac{1}{H_{m+1}} \right|\leq \frac{1}{H^2} (\ep_{J_m} +\ep_{J_{m+1}} +o(\ep_{J_m}))\leq \frac{2}{H^2} (\ep_{J_m} +o(\ep_{J_m}))$ and $\frac{1}{H_n} \leq  \frac{1}{H}(1 + \frac{\ep_{J_n}}{H}+ o(\ep_{J_n}))$, we obtain
\begin{eqnarray}
\label{eq9}
\ \ \ \  \ \ \ \ \ \ \left|\frac{\log |U_{J_n,k}|} {\log \omega_{J_n,k }(Z)} -\frac{1}{H} \right|\leq \frac{\ep_{J_n}}{H^2} +\frac{2}{H^2} \frac{ \sum_{m=0}^{n-1}( \ep_{J_m}+ o(\ep_{J_m})) \log \omega_{ J_m,K_m(k)} (Z) }{ \log \omega_{J_n,k }(Z) }
\end{eqnarray}
Let us denote $d_{m,k} = - \log \omega_{ J_m,K_m(k)} (Z)$, and $\psi_m=\ep_{J_m}$ for every $m$ and $k$. Comparing the last inequality with the desired result, all we have to show is that 
\begin{equation}
\label{conv}
\frac{ \sum_{m=0}^{n-1} \psi_{m}d_{m,k} }{ d_{n,k}} \ra 0 \ \mbox { when $n \ra +\infty$,}
\end{equation}
indepently of $k$.

This is obtained as follows: Start from $J_1$, that we suppose (without loss of generality) to be greater than 100. Recall that, by the remarks made at the beginning of Section \ref{proof}, we assumed that for every $n\geq 1$, $\eta_{J_{n }} \geq (\log J_{n }) ^{-1}$. Subsequently,   every term $J_n$ is greater than $l_n$, where $(l_n)$ is the sequence defined recursively by $l_{n+1} = l_n (1+1/\log l_n)$ and $l_1=100$. Let us study the growth rate of such a sequence. It is obvious that $\lim_{n\ra +\infty} l_n =+\infty$. We set $v_{n} = \log l_n$. We have  $v_{n+1}=v_n + \log (1+1/v_n) \geq v_n + (1-\ep)/v_n$ for every $n$, with $\ep$ that can be taken less than $1/4$ since $v_1$ is large enough. In particular, since $v_1 \geq \sqrt{2}$, $v_{2} \geq \sqrt{2} +(1-\ep)/\sqrt{2} \geq \sqrt{3}$. Recursively, if we assume that $v_n \geq \sqrt{n+1}$, then $v_{n+1} \geq \sqrt{n} +(1-\ep)/\sqrt{n} \geq \sqrt{n+1}$. Hence the sequence $(l_n)$    converges to $+\infty$ faster than $\exp{\sqrt{n}}$, and thus faster than any polynomial $n^\delta$.
(
We could be more precise, and prove using the same arguments that the growth rate of $(l_n)$ is exactly $\exp{\sqrt{n}}$.)  

Let us now find an upper bound for $\psi_m$. Using the item (1) of condition {\bf C2}, we have that 
$ \ep_j =o\left(  \frac{1}{ (\log j)^{2+\kappa} }\right)$. Using the lower bound we found for $J_m$ with $\ep$ chosen small enough, we get that $\psi_m =o(n^{-1/2(2+\kappa)}) \leq o(n^{-1-\kappa/2} ) $.

The crucial point is that $\sum_{m\geq 1} \psi_m <+\infty$. Now, since $d_{m,k}$ is a  sequence increasing  toward $+\infty$, rewrite the left term of (\ref{eq9}) as $ \sum_{m=0}^{n-1} \psi_{m}\frac{d_{m,k} }{ d_{n,k}}$, where $0\leq \frac{d_{m,k} }{ d_{n,k}}\leq 1$. By a classical Ca\!esaro method, we get that (\ref{conv}) is true, independtly of $k$.

\smallskip

This directly implies, by (\ref{eq9}), that independlty of $k$, $\left|\frac{\log |U_{J_n,k}|} {\log \omega_{J_n,k }(Z)} -\frac{1}{H} \right| \leq \kappa_n$, for some sequence $\kappa_n$ that converges to zero.

We now apply Lemma \ref{lem2}, which implies that $g$ is monofractal with exponent $H$.




\section{Around  Theorem \ref{maintheo}}
\label{secgen}

\subsection{Possible extensions for exponents greater than 1}

Let us finally say a few words about functions having regularity exponents greater than 1. The presence of a polynomial in the definition (\ref{defpoint}) of the pointwise \ho exponent is a source of problems when analyzing the local regularity after time subordination. Indeed, suppose that a continuous function $g_1$ behaves like $ |t-t_0|^{\alpha}$ ($0<\alpha<1$) around a point $t_0$, and that another continuous function  $g_2$ behaves like $a(t-t_0) + |t-t_0|^{3/2}$ ($a\neq 0$) around $t_0=g_1(t_0)$. Then $h_{g_1}(t_0) =\alpha$, $h_{g_2}(t_0)=3/2$, but $h_{g_2\circ g_1} (t_0) = \alpha$, which is different than the expected regularity $3\alpha/2$. Applying the construction above and getting a decomposition of a function $Z$ as $Z=g\circ f$, because of such problems, we didn't find any way to guarantee the monofractality of $g$.

This is related to the fact that, still for the just above toy example,   $\om_{B(t_0,r)}(g_2) \sim 2ar$ when $r$ is small enough, while one would expect $\om_{B(t_0,r)}(g_2) \sim  r^{3/2}$. The use of oscillations of   order greater than 2 (so that  $\om^2_{B(t_0,r)}(g_2) \sim r^{3/2}$) was not sufficient for us to prove Theorem \ref{maintheo} for exponents greater than 1.

An unsatisfactory result is the following: If $Z$ has all its pointwise \ho exponents less than $M>1$, then $W_{1/2M}\circ Z$ has all its exponents smaller than 1 ($W_{1/2M}$ is the Weierstrass function (\ref{defwei}) monofractal with exponent $1/2M$), and one shall try to apply Theorem \ref{maintheo} to this function.

\smallskip

As a consequence, this problem is still open and of interest.
 
\section{The case of classical monofractal functions and processes }
\label{mono}

It is satisfactory to check that classical monofractal functions verify the conditions of Theorem \ref{maintheo}, and that the exponent $H(Z)$ is actually equal to their monofractal exponent. The proofs below are also representative examples of the method used to get convergence rates for $H_j(Z)$ to $H(Z)$.

\subsection{Weierstrass-type functions}

Let $0<\alpha<1$, $\beta> \alpha$ and $b>1$ be three real numbers. 
Let $w$ be a bounded function that belongs to the global \ho class $C^\beta((0,1))$. Consider the Weierstrass-type function
\begin{equation}
\label{defwei}
Z(t) = \sum_{k=0}^{\infty} b^{-\alpha k} w(b^k t).
\end{equation}
By \cite{HB}, either the function $Z$ is $C^\beta$, or it is monofractal with exponent $\alpha$. For $w(t)=\sin(t)$, we obtain the classical Weierstrass functions monofractal with exponent $\alpha$. In fact, it is proved in \cite{HB} that, if $Z\notin C^\beta$ (which is our assumption from now on), then there is a  constant $C>1$ such that
\begin{equation}
\label{minor1}
C^{-1} 2^{-j\alpha} \leq  \om_{{j,k}} (Z) \leq C 2^{-j \alpha}.
\end{equation}
As a direct consequence,
$C^{-1/\alpha} \leq  \sum_{k=0}^{2^j-1}(\om_{j,k}(Z)) ^{1/\alpha} \leq C^{1/\alpha}$, and   obviously $H(Z) = \alpha$. 

Let us find the convergence rate of $H_j(Z)$ toward $H(Z)$.
  We are looking for a value of $\ep>0$ and for a scale $J_0$ for which $ \sum_{k=0}^{2^j-1}(\om_{j,k}(Z)) ^{1/(\alpha+\ep)}>1$, for every $j\geq J_0$. Let   $j\geq 1$. We have, by  (\ref{minor1}),  
$$ \sum_{k=0}^{2^j-1}(\om_{j,k}(Z)) ^{1/(\alpha+\ep)}  \geq 2^j  C^{-1/(\alpha+\ep)} 2^{-j \alpha/(\alpha+\ep)}  .$$
For $\ep$ small, $1/(\alpha+\ep) =1/\alpha- \ep/\alpha^2 +o(\ep)$, and thus our constraint is reached as soon as 
$1<  C^{-1/\alpha+\ep/\alpha^2+o(\ep)} 2 ^{j \ep/\alpha +o(j\ep_j)}$. This leads to
$$j \ep/\alpha +o(j\ep) + (\log_2 C)(-1/\alpha+\ep/\alpha^2+o(\ep)) >0.$$ 
There is a generation $J_0$ such that the last inequality is realized by $\ep = \frac{2   \log_2 C}{J_0}$. Subsequently, one necessarily has $1/H_j(Z) \geq 1/(\alpha +\ep)$   for every $j\geq J_0$, since 
$$\sum_{k=0}^{2^j-1}(\om_{j,k}(Z)) ^{1/(\alpha+\ep)} > \sum_{k=0}^{2^j-1}(\om_{j,k}(Z)) ^{1/H_j(Z)} =1$$
and the mapping $h \ra \sum_{k=0}^{2^j-1}(\om_{j,k}(Z)) ^{1/h}$ is increasing with $h$. Hence $H_j(Z) \leq \alpha+\ep$.

Using the same method, we obtain $H_j(Z) \geq \alpha -2\ep$ for   $j\geq J_0$.

Finally, we have found $J_0$ large enough so that for every $j\geq J_0$, $|H_j(Z)- \alpha | \leq  \ep_0 $, where we have set $\ep_{0}=2\ep=  4 \log_2 C/( J_0)$.

\medskip

For every $J\geq 1$ and $K\in\{0,...,2^J-1\}$,  we easily get the same convergence rates of $H_j(Z_{J,K})$ toward $\alpha$ from the self-affinity property of the Weierstrass functions. More precisely, fix $J$ and $K$, and let $j\geq J+1$. Remark that by construction of $Z_{J,K}$, we have $\om_{j-J,k}(Z_{J,K}) = \frac{\om_{j,K2 ^{-J}+k }(Z)}{\om_{J,K}(Z)}$. We are looking for a value of $\ep$ for which
$$ \sum_{k=0,...,2^{j}-1: I_{j,k}\subset I_{J,K}} \left(\frac{\om_{j,k}(Z)}{\om_{J,K}(Z)}\right) ^{1/(\alpha+\ep)}>1,$$
for every $j$ large enough. By  (\ref{minor1}) (used two times),  and remarking that there are $2^{j-J}$ dyadic intervals of generation $j$ included in $I_{J,K}$, we get
$$ \sum_{k=0,...,2^{j}-1: I_{j,k}\subset I_{J,K}}\left(\frac{\om_{j,k}(Z)}{\om_{J,K}(Z)}\right) ^{1/(\alpha+\ep)}  \geq 2^{j-J}  C^{-2/(\alpha+\ep)} 2^{-(j-J) \alpha/(\alpha+\ep)}  .$$
The same computations as above yield that, if we impose $\eta_J= 1/\log_2 J$ and $\ep_J= 4 (\log_2 C )(\log_2 J)/J$, then for every $j\geq J+J\eta_J$, $H_{j-J}(Z_{J,K}) \leq \alpha-\ep_J $. Similarly, we obtain $H_j(Z) \geq \alpha + 2\ep_J $  for  $j\geq J+J\eta_J$.

Finally, for every $j\geq J+J\eta_J$,   $|H_{j-J}(Z_{J,K})-\alpha | \leq  \ep_J$, and $\ep_j=o(1/(\log j)^{2+\kappa})$.
\medskip

Consequently, the Weierstrass functions satisfy {\bf C1} and {\bf C2} with $H=\alpha$, and they are also monofractal from our viewpoint.

\subsection{Sample paths of Brownian motions and fractional Brownian motions}

Classical estimations on the oscillations of sample paths of Brownian motions $(B_t)_{t\geq 0}$ yield (\cite{JAFFBEY})
\begin{eqnarray*}
\mathbb{P} \left(  \  \om_{{j,k}} (B_t) \leq j2^{-j /2}\ \right)  & \leq &  \frac{1}{2\pi} \exp {(- j^2\pi^2)}\\
\mathbb{P} \left( \ \om_{{j,k}} (B_t) \geq \frac{1}{j} 2^{-j/2} \ \right)  & \leq & \frac{4j}{2\pi} \exp {(- j^2/8)}
\end{eqnarray*}
Hence, by a classical Borel-Cantelli argument, with probability one, there is a generation $J_c$ such that for every $j\geq J_c$, we have the bounds $\frac{1}{j} 2^{-j/2}  \leq  \om_{{j,k}} (B_t) \leq j2^{-j /2}$ for the oscillations.

\medskip

The same computations as for the Weierstrass functions show that there is a generation $J_0$ such that if $j\geq J_0 \geq J_c$ , then  
\begin{eqnarray*}
\sum_{k=0}^{2^j-1} (\om_{j,k}(B_t))^{1/(1/2+\ep_{J_0})} >1 \mbox{ and } \sum_{k=0}^{2^j-1} (\om_{j,k}(B_t))^{1/(1/2-\ep_{J_0})}  <1,
\end{eqnarray*}
where $\ep_{J_0} \geq C \frac{\log J_0}{J_0}$ (for some suitable constant $C$). As a consequence, $|H_j(B_t) -1/2| \leq \ep_{J_0}$ for every $j\geq J_0$.

\medskip

The self-similarity property of Brownian motions yields that for every $J\geq 1$ and $K\in\{0,...,2^J-1\}$,   for every $j\geq J/\log J$, $|H_j((B_t)_{J,K})-1/2| \leq \ep_J $, where $\ep_J =C   \frac{\log^2 J}{J}$, for some constant $C $ independent of $J$ and $K$. We omit the details here, that can be easily checked by the reader.

\medskip

 Consequently, a sample path of Brownian motion satisfies with probability one {\bf C1} and {\bf C2}, with $H(B_t)=1/2$.

\medskip

Similar estimations on the oscillations of fractional Brownian motions $B_h$ of Hurst exponent $h$ lead to the same almost sure result for the sample paths, which also satisfy  almost surely  {\bf C1} and {\bf C2} with $H(B_h) =h$.

%


\section{Applications to self-similar functions: Theorem \ref{thself}}
\label{self}
 
 We consider the class of self-similar functions defined in Definition \ref{defiself}, with the parameters of $Z$ and the contractions $S_k$ satisfying  (\ref{cond1}). 
 
%




 The multifractal analysis of such a function $Z$ is performed in \cite{JAFFFORM}. Here we are going to prove that, under the  conditions    (\ref{cond1}) on the $\la_k$ and the $S_k$, $Z$ is CMT, and that the multifractal behavior of $Z$ can be directly deduced from this analysis. It is a case where our analysis provides a   natural way to compute the singularity spectrum  of $Z$.

\subsection{Preliminary results on the oscillations of $Z$}

Let us introduce some notations: for every $n\geq 1$, for every $(\ep_1,\ep_2,...,\ep_n) \in \{0,1,...,d-1\} ^n$, we denote $I_{\ep_1,\ep_2,...,\ep_n}$ the interval $S_{\ep_1}  \circ S_{\ep_2}  \circ... \circ  S_{\ep_n}  (\zu)$. The integer $n$ being given, the open intervals $\stackrel{\circ}{(I_{\ep_1,\ep_2,...,\ep_n})}$ are pairwise disjoint, and the union of the closed intervals $ I_{\ep_1,\ep_2,...,\ep_n}$  equals $\zu$. Now fix an integer $n\geq 1$ and a sequence $(\ep_1,\ep_2,...,\ep_n) \in \{0,1,...,d-1\} ^n$.  The interval  $ I_{\ep_1,\ep_2,...,\ep_n}$ has a length equal to $r_{\ep_1}r_{\ep_2}\cdot \cdot\cdot r_{\ep_n}$.  
Finally, by iterating $n$ times formula (\ref{defself}), we get that for every $t\in  I_{\ep_1,\ep_2,...,\ep_n}$,
\begin{eqnarray}
\label{decompf}   Z(t) & = &\la_{\ep_1}\cdot\la_{\ep_2}\cdot\cdot\cdot\la_{\ep_n} \cdot (Z\circ S_{\ep_n}^{-1} \circ S_{\ep_{n-1}}^{-1} \circ... \circ  S_{\ep_1}^{-1} ) (t) \\
 \nonumber&+ &\la_{\ep_1}\cdot\la_{\ep_2}\cdot\cdot\cdot\la_{\ep_{n-1}} \cdot (\phi\circ S_{\ep_{n-1}}^{-1} \circ S_{\ep_{n-2}}^{-1} \circ... \circ  S_{\ep_{1}}^{-1} ) (t) \\ 
 \nonumber &+& ... \\
 \nonumber &+&  \la_{\ep_1} \cdot\la_{\ep_2}\cdot (\phi\circ S_{\ep_2}^{-1} \circ S_{\ep_1}^{-1}  ) (t) \\
 \nonumber &+&  \la_{\ep_1} \cdot (\phi\circ S_{\ep_1}^{-1}  ) (t) \\
 \nonumber &+&  \phi(t).
 \end{eqnarray}


Recall that $\chi_{\max}$ is defined in (\ref{cond1}).
\begin{proposition}
\label{lem3}
Let $\kappa =   \frac{\chi_{\max}}{1-\chi_{\max}}$. Then either $Z$ is a $\kappa$-Lipschitz function, or there is a constant $C>1$ such that for every $n\geq 1$, for every $(\ep_1,\ep_2,...,\ep_n) \in \{0,1,...,d-1\} ^n$, 
\begin{equation}
\label{majmin1}
C^{-1}  \cdot  |\la_{\ep_1}\cdot\la_{\ep_2}\cdot\cdot\cdot\la_{\ep_n}  |\leq \om_{I_{\ep_1,\ep_2,...,\ep_n}}(Z) \leq C \cdot |\la_{\ep_1}\cdot\la_{\ep_2}\cdot\cdot\cdot\la_{\ep_n} |.
\end{equation}
\end{proposition}

 \begin{proof}
We first find an upper-bound for $ \om_{I_{\ep_1,\ep_2,...,\ep_n}}(Z) $. We use the iterated formula (\ref{decompf}).

Let $n$ and $(\ep_1,\ep_2,...,\ep_n) \in \{0,1,...,d-1\} ^n$. 
Remark that when $t$ ranges in ${I_{\ep_1,\ep_2,...,\ep_n}}$,   $(S_{\ep_n}^{-1} \circ S_{\ep_{n-1}}^{-1} \circ... \circ  S_{\ep_1}^{-1} ) (t) $ ranges in $\zu$. Hence the oscillation of the first term of (\ref{decompf}) is upper-bounded by $ | \la_{\ep_1} \la_{\ep_2 } \!  \cdot  \! \cdot \! \cdot \! \la_{\ep_n}  | \cdot \om_{\zu}(Z)$.

Now, for every $k\in\{1,...,n-1\}$, when $t$ ranges in ${I_{\ep_1,\ep_2,...,\ep_n}}$, $(S_{\ep_k}^{-1} \circ S_{\ep_{k-1}}^{-1} \circ... \circ  S_{\ep_1}^{-1} ) (t)  $ ranges in $I_{\ep_{k+1},...,\ep_n}$. Using that $\phi$ is  a  Lipschitz function, we get that the oscillation of each term of the form  $ \la_{\ep_1}\la_{\ep_2} \!  \cdot  \! \cdot \! \cdot \! \la_{\ep_{k}} \cdot (\phi \circ S_{\ep_{k}}^{-1} \circ S_{\ep_{k-1}}^{-1} \circ... \circ  S_{\ep_{1}}^{-1} ) (t) $ is  upper bounded by $ |\la_{\ep_1} \la_{\ep_2}   \!  \cdot  \! \cdot \! \cdot \!  \la_{\ep_{k}} |  (    r_{\ep_{k+1}}  \!  \cdot  \! \cdot \! \cdot \! r_{\ep_n})$.
Finally, we obtain using (\ref{cond1})
\begin{eqnarray*}
\om_{I_{\ep_1,\ep_2,...,\ep_n}}(Z)  & \leq  &  | \la_{\ep_1} \la_{\ep_2}  \!  \cdot  \! \cdot \! \cdot \!  \la_{\ep_n}  |+   \sum_{k=1} ^{n-1}  |\la_{\ep_1} \la_{\ep_2} \!  \cdot  \! \cdot \! \cdot \!  \la_{\ep_{k}}  |\cdot (   r_{\ep_{k+1}}  \!  \cdot  \! \cdot \! \cdot \! r_{\ep_n})\\
&\leq &   | \la_{\ep_1} \la_{\ep_2}  \!  \cdot  \! \cdot \! \cdot \!  \la_{\ep_n}  |\Big[ 1+  \sum_{k=1}^{n-1}\big( \prod_{j=1}^{k} \frac{r_j}{\la_j} \Big) \Big] \\
&\leq &  |  \la_{\ep_1} \la_{\ep_2}  \!  \cdot  \! \cdot \! \cdot \! \la_{\ep_n} | \Big[ 1+  \sum_{k=1}^{n-1} \chi_{\max}^k \Big] \leq C_1 |  \la_{\ep_1}\la_{\ep_2}\ \!  \cdot  \! \cdot \! \cdot \!  \la_{\ep_n}  |,
\end{eqnarray*} 
where $C_1= 1+  \sum_{k=1}^{+\infty} \chi_{\max}^k
< +\infty$.

\smallskip

We now move to the lower bound. Assume that $Z$ is not $\kappa$-Lipschitz. There are two real numbers $0\leq t_0 ,t'_0 \leq 1$ such that $|Z(t'_0)-Z(t_0)| \geq (\kappa+\eta) |t'_0-t_0|$, for some $\eta>0$. Let $n$ and $(\ep_1,\ep_2,...,\ep_n) \in \{0,1,...,d-1\} ^n$. 
Let us call $t_n = S_1 \circ S_{2} \circ ... \circ S_{n}(t_0)$ and  $t'_n = S_1 \circ S_{2} \circ ... \circ S_{n} (y_0)$.
We obviously have $t_n, t'_n \in I_{\ep_1,\ep_2,...,\ep_n}$, and thus $\om_{I_{\ep_1,\ep_2,...,\ep_n}}(Z)  \geq |Z(t'_n) -Z(t_n)|$. Using again   (\ref{decompf}), we get by the same lines of computations as above
\begin{eqnarray*}
&&\hspace{-10mm} |Z(t'_n) -Z(t_n)|  \\ 
\geq  && \hspace{-5mm}    |\la_{\ep_1} \la_{\ep_2}  \!  \cdot  \! \cdot \! \cdot \!  \la_{\ep_n}  |\cdot |Z(t'_0) -Z(t_0)|    -   \sum_{k=1} ^{n-1}  |\la_{\ep_1} \la_{\ep_2} \!  \cdot  \! \cdot \! \cdot \!  \la_{\ep_{k}}  | (   r_{\ep_{k+1}}  \!  \cdot  \! \cdot \! \cdot \! r_{\ep_n}) |t'_0-t_0|\\
\geq &&\hspace{-5mm}  |   \la_{\ep_1} \la_{\ep_2}  \!  \cdot  \! \cdot \! \cdot \!  \la_{\ep_n} |\cdot  |Z(t'_0) -Z(t_0)| \left[ 1 -  \sum_{k=1}^{n-1}\big( \prod_{j=1}^{k} \frac{r_j}{ |\la_j |} \Big)  \frac{|t'_0 -t_0|}{|Z(t'_0) -Z(t_0)|}\right] \\
\geq &&\hspace{-5mm}    |  \la_{\ep_1} \la_{\ep_2}  \!  \cdot  \! \cdot \! \cdot \! \la_{\ep_n}  |\cdot |Z(t'_0) -Z(t_0)|  \Big[ 1- \frac{ 1} {\kappa +\eta}\sum_{k=1}^{+\infty} \chi_{\max}^k \Big] \geq C_2 \cdot  | \la_{\ep_1}\la_{\ep_2}\ \!  \cdot  \! \cdot \! \cdot \!  \la_{\ep_n}   |,\end{eqnarray*} 
where $C_2 = |Z(t'_0) -Z(t_0)|  ( 1- \frac{ 1} {\kappa +\eta}\frac{\chi_{\max}}{1-\chi_{\max}})  >0$ by assumption.

Finally, (\ref{majmin1}) is proved with $C= \max(C_1, C_2^{-1})$.
\end{proof}

\subsection{Comparaison of $Z$ with a self-similar measure}

In order to prove that the function $Z$ (\ref{defself}) satisfies our conditions  {\bf C1} and {\bf C2}, we introduce a self-similar measure $\mu$, whose multifractal behavior will be compared with the one of  $Z$, and the notion of multifractal formalism.

Let us consider the  exponent $\beta>1$ such that (\ref{defbeta}) holds and the associated self-similar measure $\mu$ defined by (\ref{defmu}) $  \mu = \sum_{k=0} ^{d-1} p_k \cdot  (\mu\circ S_k ^{-1}).$ In our case where the similitudes do not overlap, it  is easily checked that by construction, for every $n$ and $(\ep_1,\ep_2,...,\ep_n) \in \{0,1,...,d-1\} ^n$, we have $\mu(  I_{\ep_1,\ep_2,...,\ep_n}) = p_{\ep_1} p_{\ep_2}  \!  \cdot  \! \cdot \! \cdot   p_{\ep_n} =  | \la_{\ep_1} \la_{\ep_2}  \!  \cdot  \! \cdot \! \cdot \!  \la_{\ep_n} |^{\beta}$.

This class of measures has been extensively studied  \cite{BMP,CM,OLSEN,PERES}. For instance, the multifractal analysis of $\mu$ is very well known. For this, let us introduce the so-called $L^q$-spectrum of $\mu$ defined by
\begin{equation}
\label{deftau}
\tau_\mu: q\in \R \mapsto  \tau_\mu(q)= \liminf_{j\ra +\infty}   \tau_\mu(j,q), \ \mbox{ where } \tau_\mu(j,q) =  \frac {\log_2 \sum_{k=0} ^{2^j-1} \mu(I_{j,k}) ^q}{-  j}.
\end{equation}
We only recall the properties we need \cite{CM,OLSEN,PERES}
 \begin{proposition}
 \label{prop4}
  \begin{enumerate}
 \item
 For every $q\in \R$, $ \tau_\mu(q)$  is the unique real number satisfying the equation $\sum_{k=0}^{d-1} (p_k )^q (r_k)^{\tau_\mu(q)} =1$.  The mapping $q\mapsto \tau_\mu(q)$  is analytic on its support. Moreover, the liminf used to define $\tau_\mu(q)$ is in fact a limit for every $q$ such that $\tau_\mu(q)$ is finite. 

\item
There is an interval of exponents $I_\mu= [\alpha_{\min},\alpha_{\max}]$ such that for every $\alpha\in I_\mu$, $\tilde d_{\mu}(\alpha) = (\tau_\mu)^*(\alpha) $, where $(\tau_\mu)^*(\alpha) := \inf_{q\in \R} (q\alpha- \tau_\mu(q))$ is by definition the Legendre transform of $\tau_\mu$. 

\item
If $\alpha\notin I_\mu$, then $ \{x : \alpha_\mu(t)=\alpha\} = \emptyset $.

\item
There is  $M \geq 1$ such that for every $j,k$ large enough, $2^{-j M} \leq \mu(I_{j,k}) \leq 2 ^{-j /M}$.
\end{enumerate}
\end{proposition}
Part (2) above is known as the multifractal formalism for measures, when it holds.
 
 \smallskip
 
 Let us come back to the function $Z$. The reader can check that such a function $Z$ satisfies {\bf C1} and {\bf C2}, and is thus CMT. Here we propose a quick proof of Theorem \ref{thself}, especially adapted to this case.
 
 The aim is to prove that $Z$ can be written $Z= g\circ$ 
 Each dyadic interval $I_{j,k}$ is included in one dyadic  interval $I_{\ep_1,...,\ep_n}$, and contains a dyadic interval $I_{\ep_1,...,\ep_n,\ep_{n+1}}$, such that $I_{\ep_1,...,\ep_n}$ and $I_{\ep_1,...,\ep_n,\ep_{n+1}}$ can be written respectively $I_{j',k'}$ and $I_{j'',k''}$ with $0\leq j-j' , j''-j'  \leq C$, for some constant $C$ independent of $j$ and $k$. Consequently, $\om_{ I_{\ep_1,...,\ep_n,\ep_{n+1}}} (Z) \leq \om_{I_{j,k} }  (Z) \leq \om_{ I_{\ep_1,...,\ep_n}}  (Z) $, and thus
 $$ C^{-1}\mu({ I_{\ep_1,...,\ep_n,\ep_{n+1}}})^{1/\beta} \leq \om_{{j,k} }  (Z) \leq C\mu({ I_{\ep_1,...,\ep_n}} )^{1/\beta}.
 $$ 
 
 Using now   the self-similarity properties of the measure and the open set condition,  we see that   $ \max_k ({p_k}) \cdot {\mu (I_{\ep_1,...,\ep_n,\ep_{n+1}} ) }  \geq  {\mu(I_{j,k}) }$ and $\mu(I_{j,k}) \geq  \min_k (p_k) \cdot \mu({ I_{\ep_1,...,\ep_{n}}})$. Hence, combining this with the last double inequality, we obtain  that for every $j$ and $k$   
 \begin{equation}
 \label{majmin3}
  C^{-1}\mu(I_{j,k})^{1/\beta} \leq \om_{{j,k} }  (Z) \leq C\mu( I_{j,k})^{1/\beta} 
  \end{equation}
  for another constant $C$, i.e. Proposition \ref{lem3} extends to all dyadic intervals.

\smallskip

Let us check that $Z$ satisfies conditions {\bf C1} and {\bf C2}. 
Remark that, because of (\ref{majmin3}),
$$ C^{-\beta} =  C^{-\beta}  \sum_{k=0}^{2^j-1} \mu(I_{j,k}) \leq \sum_{k=0}^{2^j-1} (\om_{j,k}(Z))^{\beta} \leq C^{\beta}  \sum_{k=0}^{2^j-1} \mu(I_{j,k}) = C^{\beta}  .$$ 
Let $\ep_1>0$. Using part (4) of Proposition \ref{prop4} to find upper- and lower-bounds for $\om_{j,k}(Z)$ uniformly in $k$, we obtain
\begin{eqnarray*}
&&  \sum_{k=0}^{2^j-1} (\om_{j,k}(Z))^{\beta -\ep_1} \geq    \sum_{k=0}^{2^j-1} (\om_{j,k}(Z))^{\beta} (C\mu(I_{j,k}))^{ -\ep_1} \geq C^{-\beta-\ep_1}  2^{j\ep_1/M} \\
 \mbox{and}&& \sum_{k=0}^{2^j-1} (\om_{j,k}(Z))^{\beta +\ep_1} \leq    \sum_{k=0}^{2^j-1} (\om_{j,k}(Z))^{\beta} (C\mu(I_{j,k}))^{ \ep_1} \geq C^{\beta+\ep_1}  2^{-j\ep_1 M}  .
\end{eqnarray*}
Hence, the same computations as in the case of Weierstrass functions lead to the following choice: for some $J_0$ large enough,   we set $\ep_{1} = \frac{2\beta M \log_2 C}{\log J_0}$, and thus for every $j \geq J_0$, $|H_j(Z) - H(Z)|  \leq \ep_1$.

Let now $J,K$ be two integers, $\ep >0$, and focus on $H(Z_{J,K})$. The same computations as above  and as in the Weierstrass case yield for $j\geq J$
\begin{eqnarray*}
\sum_{k'=0}^{2^{j-J}-1} ({\om_{j-J,k'}(Z_{J,K})})^{\beta -\ep} &= &  \sum_{k=0: I_{j,k}\subset I_{J,K}}^{2^j-1} \left(\frac{\om_{j,k}(Z)}{\om_{J,K}(Z)}\right)^{\beta -\ep} .
\end{eqnarray*}
First notice that $\left(\frac{1}{\om_{J,K}(Z)}\right)^{\beta -\ep} \geq   \left(\frac{1}{C\mu(I_{J,K})^{1/\beta}}\right)^{ -\beta+\ep} \geq C^{-\beta-\ep}   \mu(I_{J,K}) ^{-1+\ep/\beta}
$.
Then we remark that 
\begin{eqnarray*}
 \sum_{k=0: I_{j,k}\subset I_{J,K}}^{2^j-1} (\om_{j,k}(Z))^{\beta -\ep}  & \geq &\sum_{k=0: I_{j,k}\subset I_{J,K}}^{2^j-1} (\om_{j,k}(Z))^{\beta } (C\mu(I_{j,k}))^{ -\ep} .
 \end{eqnarray*}
Combining these inequalities we get
\begin{eqnarray*}
\sum_{k'=0}^{2^{j-J}-1} ({\om_{j-J,k'}(Z_{J,K})})^{\beta -\ep}&\geq &  C^{-\beta-2\ep} \sum_{k=0: I_{j,k}\subset I_{J,K}}^{2^j-1} \frac{(\om_{j,k}(Z))^{\beta }}{\mu(I_{J,K})}   \frac {\mu(I_{J,K}) ^{\ep/\beta} } {    \mu(I_{j,k})^{ \ep}   }.\end{eqnarray*}
Let us focus on $\frac {\mu(I_{J,K}) ^{\ep/\beta} } {    \mu(I_{j,k})^{ \ep}   }$. Since $\beta>1$, we have 
$\frac {\mu(I_{J,K}) ^{\ep/\beta} } {    \mu(I_{j,k})^{ \ep}   } \geq \left(\frac {\mu(I_{J,K})  } {    \mu(I_{j,k})}\right)^{ \ep}   $. This quantity is lower bounded by $L^{(j-J) \ep}$ for some constant $L$ (uniformly in $k$ and $K$), since the ratio of the $\mu$-measures of a dyadic interval and its father (in the dyadic tree) is uniformly upper- and lower-bounded for our dyadic self-similar measure $\mu$. Finally, we obtain
\begin{eqnarray*}
\sum_{k'=0}^{2^{j-J}-1} ({\om_{j-J,k'}(Z_{J,K})})^{\beta -\ep}&\geq &  C^{-\beta-2\ep} \sum_{k=0: I_{j,k}\subset I_{J,K}}^{2^j-1} \frac{(\om_{j,k}(Z))^{\beta }}{\mu(I_{J,K})} L^{(j-J) \ep}\\
&\geq &  C^{-2\beta-2\ep} \sum_{k=0: I_{j,k}\subset I_{J,K}}^{2^j-1} \frac{\mu(I_{j,k})}{\mu(I_{J,K})} L^{(j-J) \ep} \\
& \geq & C^{-2\beta-2\ep} L^{(j-J) \ep}.\end{eqnarray*}
Hence if we fix $\eta_J = 1/\log_2 J$, then the  sum above is greater than 1 as soon as $j\geq J+[J\eta_J]$ and $\ep \geq  \frac{4 \beta \log C \log_2 J}{J \log L}$. Thus  for $j\geq J+[J\eta_J]$, $H_{j-J}(Z_{J,K}) - H(Z_{J,K}) \leq \frac{1}{\beta-\ep} \leq \frac{1}{\beta} +\ep_J$ with $\ep_J=  \frac{8 \log C \log_2 J}{J  \beta\log L}$. 

Similarly one shows that  for $j\geq J+[J\eta_J]$, $H_{j-J}(Z_{J,K}) - H(Z_{J,K}) \geq \frac{1}{\beta} -\ep_J$,  and {\bf C2} holds true for $Z$.

\medskip

Applying Theorem \ref{maintheo} yields that $Z$ is CMT and can be written as  $Z=g\circ F$, where $g$ is monofractal of exponent $1/\beta$.

\subsection{Computation of the singularity spectrum of $F$}

Applying directly the construction of Section \ref{proof}, we find a function $g$ monofractal with exponent $1/\beta$ and a strictly increasing function $f$ such that $Z= g \circ f$.

One can even enhance this result as follows. Following the proof of Section \ref{proof}, we see that for every $p\geq 1$, for every $k\in \{0,1,...,2^{J_p}-1\}$, $$ \mu(I_{J_p,k}) ^{1 + \kappa_p} \leq |f(I_{J_p,k}) | \leq  \mu(I_{J_p,k}) ^{1 -\kappa_p},$$
where $(\kappa_p)_{p\geq 1}$ is a positive sequence decreasing to zero and $\mu$ is defined by (\ref{defmu}).

Let us denote by $F$ the integral of the self-similar measure $\mu$, i.e. for $t\in\zu$ $F(t)= \mu([0,x])$. We claim that $f= g_1 \circ F$ for some function $g$ which belongs to $C^{1-\eta}(\zu)$, for every $\eta>0$. Indeed,  define for every $t\in\zu$ $g_1(t) = f\circ F^{-1} (t)$. This is possible since $F$ is an homeomorphism of $\zu$. 

By construction, for every $p\geq 1$,   for  every $k\in \{0,1,...,2^{J_p}-1\}$, $g_1(F \big( I_{J_p,k}) \big) = f \circ F^{-1} \circ  F \big( I_{J_p,k})    = f  \big( I_{J_p,k})  $, thus by the inequality above, 
$$ [F(I_{J_p,k} )|^{1 + \kappa_p} =  \mu(I_{J_p,k}) ^{1 + \kappa_p} \leq |g_1(F  \big( I_{J_p,k}) \big)  | \leq  \mu(I_{J_p,k}) ^{1 -\kappa_p} = [F(I_{J_p,k} )|^{1 -\kappa_p} ,$$
where we used that $[F(I_{J_p,k} )|  =  \mu(I_{J_p,k})$. Now the sets of intervals $\{ F \big( I_{J_p,k}) \big): k\in \{0,1,...,2^{J_p}-1\}\}$ obviously forms a covering of $\zu$ to which Lemma \ref{lem2} can be applied with $H=1$.

FInally, we find that $Z=g\circ g_1 \circ F = g_2 \circ F$, where $g_2$ is clearly monofractal with exponent $1/\beta$ since $g$ and $g_1$ are  monofractal respectively with exponents $1/\beta$ and $1$ (the resulting function $g_2$ is monofractal since the oscillations of $g$ and $g_1$ are upper and most important lower bounded on every interval). 


\section{An example of function satisfying {\bf C1-C2} in a triadic basis}
\label{multi}

We recall the contruction of multifractal functions of \cite{OKA}, which somehow generalizes the Bourbaki's and Perkin's functions. 

\begin{center}
\begin{figure}
\label{fig1}
\includegraphics[width=10cm,height=5cm]{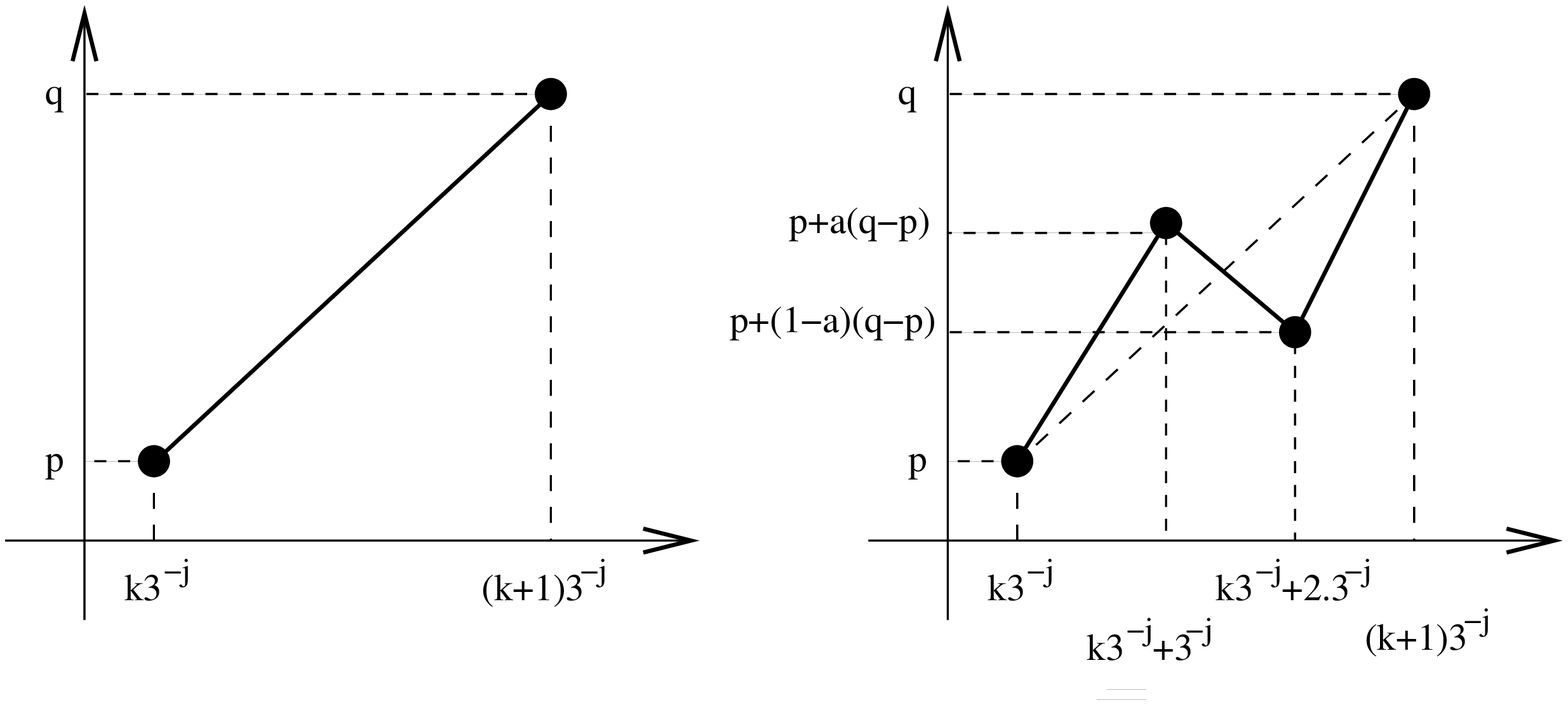}
\caption{Iterated construction of $Z_a$, from step $j$ to step $j+1$}
\end{figure}
\end{center}

Let us consider the function $Z_a$ defined for $0\leq a\leq 1$ as the limit of an iterated construction: Start from $Z_a^0(t)=t $ on $\zu$, and define $Z_a^j(t)$ recursively on $\zu$ by the following scheme: Suppose that $Z_a^j$ is continuous and piecewiese affine on each triadic interval $[k3^{-j}, (k+1)3^{-j}]$, $k\in\{0,...,3^j-1\}$. Then  $Z_a^{j+1}$ is constructed as follows: On each triadic interval $[k3^{-j}, (k+1)3^{-j}]$, $Z_a^{j+1}$ is still a continuous function which is affine on each triadic subinterval  $[k'3^{-(j+1)}, (k'+1)3^{-(j+1)}]$ included in $[k3^{-j}, (k+1)3^{-j}]$,   and 
\begin{eqnarray*}
Z_a^{j+1} (k3^{-j}) & =  & Z_a^{j} (k3^{-j})\\
Z_a^{j+1} (k3^{-j}+ 3^{-(j+1)}) & =  & Z_a^{j} (k3^{-j}) + a \Big( Z_a^{j} ((k+1)3^{-j})- Z_a^{j} (k3^{-j}) \Big) \\
Z_a^{j+1} (k3^{-j}+ 2.3^{-(j+1)}) & =  & Z_a^{j} (k3^{-j})+ (1-a) \Big( Z_a^{j} ((k+1)3^{-j})- Z_a^{j} (k3^{-j}) \Big) \\
Z_a^{j+1} ((k+1)3^{-j}) & =  & Z_a^{j} ((k+1)3^{-j}).
\end{eqnarray*}
This simple construction is better explained by the Figure \ref{fig1}.

It is straightforward to see that the sequence $(Z_a^j)_{j\geq 1}$ converges uniformly to a continuous function $Z_a$ as soon as $0<a<1$.
Bourbaki's function is obtained when $a=2/3$, while Perkin's function corresponds to $a=5/6$.

For $a\leq 1/2$, the function is simply the integral of a trinomial measure of parameters $(a, 1-2a,a)$, hence its singularity spectrum is completely known. We are going to explain why the functions $Z_a$, when $a\geq 1/2$, satisfy our assumptions, and thus can be written as the composition of a monofractal function $g$ (with an exponent $H$ we are going to determine) with an increasing function. We will also deduce from this study the singularity spectrum of $Z_a$.

\medskip

For $a>1/2$, the limit function $Z_a$ is nowhere monotone. 
Let us compute the oscillations of $Z_a$ on each triadic interval. 

Remark first that  the slope of $Z_a^1$ on $[0,1/3]$ is $3a$, it is $-3(2a-1)$ on $[1/3,2/3]$ and $3a$ on $[2/3,1]$. Iteratively, if $j \geq 1$ and $k\in\{0,...,3^j-1\}$, we write $k3^{-j} = \sum_{p=1}^j \xi_p 3^{-p}$, with $\xi_i \in\{0,1,2\}$. Then the slope of $Z_a^j$ on   $[k3^{-j}, (k+1)3^{-j}]$ is simply
$$(3a) ^{n_{k,j,0}} (-3(2a-1))^{n_{k,j,1}}(3a) ^{n_{k,j,2}} = 3^j (a) ^{n_{k,j,0}} (-(2a-1))^{n_{k,j,1}}(a) ^{n_{k,j,2}},$$
where $n_{k,j,i}$ is the number of integers $p\in\{1,...,j\}$ such that $\xi_p=i$ (for $i=0,1,2$) in the triadic decomposition of $k3^{-j}$.
 
 Let us consider the trinomial measure $\mu_a$ of parameters $(\frac{a}{4a-1},\frac{2a-1}{4a-1},\frac{a}{4a-1})$. Then it is obvious that the absolute value of the slope of $Z_a^j$ on each triadic interval $[k3^{-j}, (k+1)3^{-j}]$ can be written as
 $ \mu_a([k3^{-j}, (k+1)3^{-j}])3^j (4a-1)^j.$
 As a final remark, we also notice that the oscillations of $Z_a$ on each triadic interval $[k3^{-j}, (k+1)3^{-j}]$ is the same as the oscillations of $Z^j_a$ on each triadic interval $[k3^{-j}, (k+1)3^{-j}]$, which is equal to $3^{-j}$ times the slope, i.e. 
 \begin{equation}
 \label{eq11}
  \mu_a([k3^{-j}, (k+1)3^{-j}])  (4a-1)^j.
  \end{equation}

Let $q\in\R$. Let us compute the sum of the oscillations of $Z_a$ at generation $j$. We have
\begin{equation}
\label{eq10}
\sum_{k=0}^{3^j-1} (\om_{[k3^{-j}, (k+1)3^{-j}]}(Z_a))^{q} =\sum_{k=0}^{3^j-1}  (\mu_a([k3^{-j}, (k+1)3^{-j}]))^{q} 3^{qj\log_3  (4a-1)}.
\end{equation}
Let us explain now how we easily compute  the exponent $H_a$ such that (\ref{defh}) holds true.
For a multinomial measure $\mu_a$ (in fact, for any positive Borel measure), it is very classical in multifractal analysis to introduce the functions $\tau_{\mu_a,j}(q)$ and the scaling function $\tau_{\mu_a}(q)$ defined for $q\in\R$ as (\ref{deftau}) but in the triadic basis:
\begin{eqnarray*}
\tau_{\mu_a}(q)= \liminf_{j\ra +\infty} \tau_{\mu_a,j}(q), \ \mbox{ where } \ \tau_{\mu_a,j}(q)=  \frac{\log_3 \sum_{k=0}^{3^j-1}  (\mu_a([k3^{-j}, (k+1)3^{-j}]))^{q}}{ -j}.
\end{eqnarray*}
In our simple case, it is easy to see that for every $j\geq 1$ and $q\in\R$
\begin{eqnarray*} \tau_{\mu_a,j}(q) \ =  \ \tau_{\mu_a}(q) & = & -\log_3\left(\left(\frac{a}{4a-1}\right)^q+\left( \frac{2a-1}{4a-1}\right)^q+\left(\frac{a}{4a-1}\right)^q\right)\\
 & =  & -\log_3(2(a)^q+(2a-1)^q) + q\log_3 (4a-1).
 \end{eqnarray*}
 
What matters to us is the value of $q$ for which the sum in (\ref{eq10}) equals 1. Let us write this specific value $q$ as $1/H$, for some $H>0$. When this sum is 1, then we have
$$ 3^{-j\tau_{\mu}(1/H)}  3^{ j(\log_3  (4a-1))/H} = \sum_{k=0}^{3^j-1}  (\mu_a([k3^{-j}, (k+1)3^{-j}]))^{1/H} 3^{j(\log_3  (4a-1))/H}=1.$$
Let $H_a$ be the   solution of the equation $-\tau_{\mu}(1/H_a) +  \log_3  (4a-1)/H_a=0$, which is equivalent to 
\begin{equation}
\label{defha}
 2(a)^{1/H_a}+(2a-1)^{1/H_a} =1.
 \end{equation}
  This solution is positive, unique, and strictly smaller than 1. Hence, in this case, the monofractal exponent $H_a$ is defined through an implicit formula.


In order to get the whole condition {\bf C2}, it suffices to notice that any rescaled function $(Z_a)_{J,K}$ (as defined in (\ref{eq00}), but here with triadic intervals) is actually equal to $Z_a$ (if $Z_a$ is increasing on $[k3^{-j},(k+1)3^{-j}]$) or to $Z_a(1-.)$ (if $Z_a$ is decreasing on $[k3^{-j},(k+1)3^{-j}]$). Hence $H((Z_a)_{J,K})$ is a limit for every $J,K$, and is even constant equal to $H_a$. Thus {\bf C2} is satisfied.

\medskip

\begin{center}
\begin{figure}
\begin{center}
\label{fig2}
\hspace{-15mm} \includegraphics[width=14cm,height=5cm]{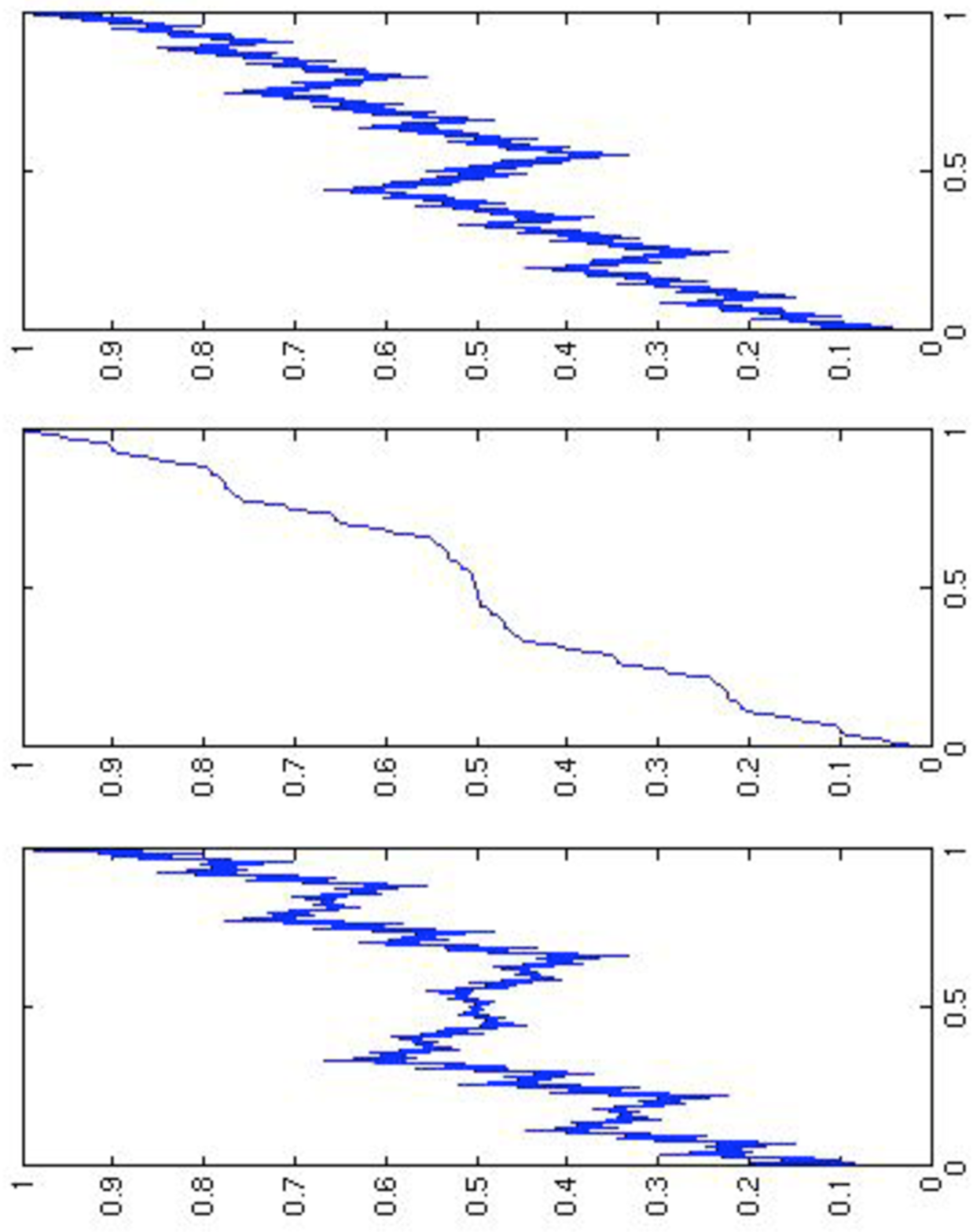} \\
\hspace{-15mm} \includegraphics[width=12cm,height=4cm]{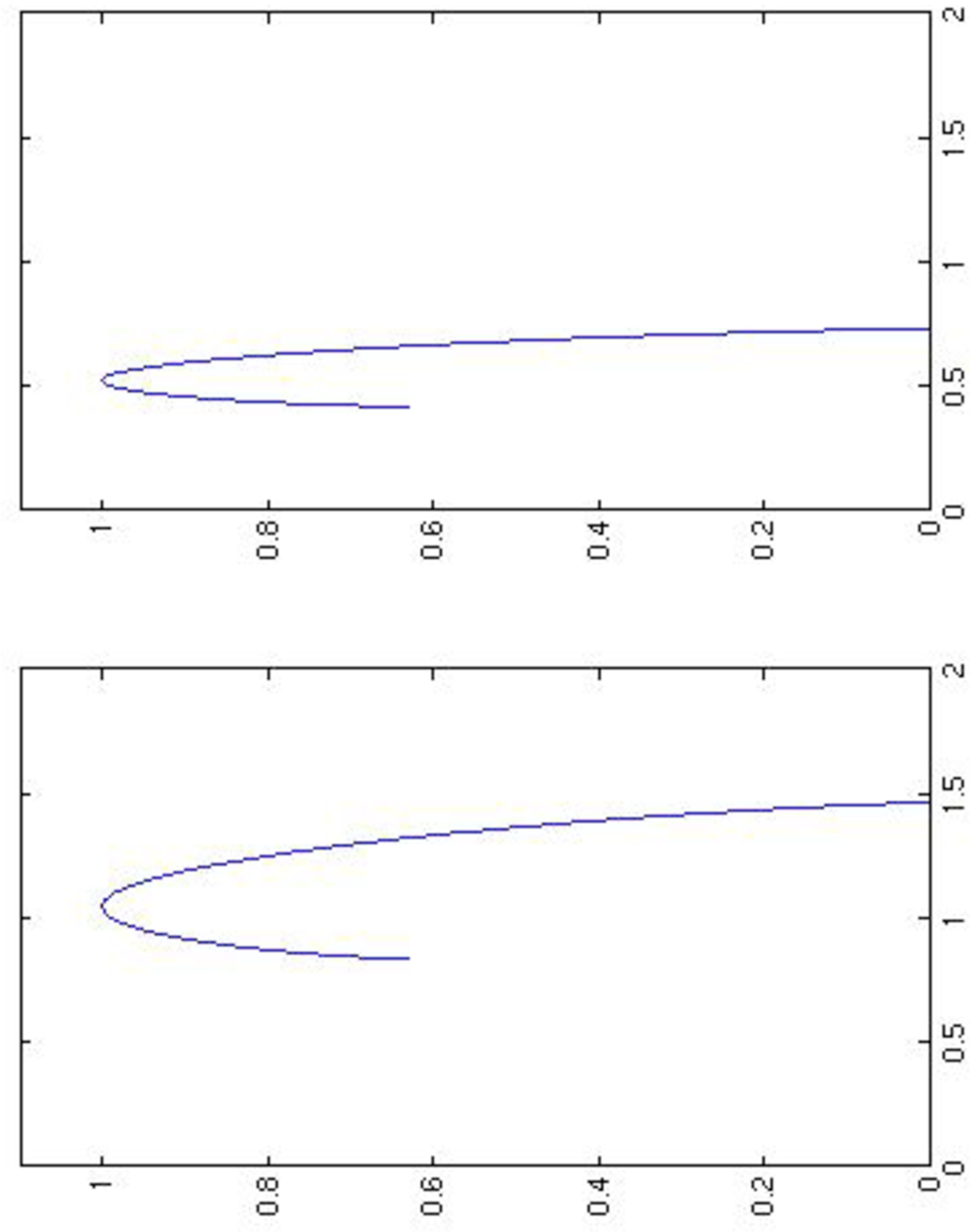} 
\caption{{\bf Top:} Bourbaki's function $f_{2/3}$ on the left, the multifractal time change $F_{2/3}$ in the middle, and on the right the monofractal function $g_{2/3}$ of exponent $1/2$ such that $f_{2/3} = g_{2/3}\circ F_{2/3}$. {\bf  Bottom:} Singularity spectra of $\mu_{2/3}$ on the left, of $f_{2/3}$ on the right.}
\end{center}
\end{figure}
\end{center}

We can then apply Theorem \ref{maintheo}, and $Z_a$ is the composition of a monofractal function $g_a$ of exponent $H_a$ with an increasing function $Z_a$. For $a=2/3$, we see that $H_a=1/2$ is the solution to (\ref{defha}), since $2(2/3)^{2} + (1/3)^{2} =1$. We have plotted in Figure \ref{fig2} the Bourbaki's function $f_{2/3}$, its corresponding time change $F_{2/3}$ and the corresponding monofractal function $g_{2/3}$ of exponent $1/2$ such that $f_{2/3}=g_{2/3}\circ F_{2/3}$.

\medskip
 
 In this case, we can even go further and compute the singularity spectrum of $Z_a$.  
The trinomial measure satisfy the multifractal formalism for measures, i.e. the singularity spectrum of $\mu$ is given by the Legendre transform of $\tau_{\mu_a}$:
 $$ d_{\mu_a}(\alpha) = (\tau_{\mu_a})^*(\alpha):= \inf_{q\in\R} (q\alpha-\tau_{\mu_a}(q)),$$
for every $\alpha \in [-\log_3 (a/(4a-1)), -\log_3(2a-1)/(4a-1)]$.

%

It is easy to see, using (\ref{eq11}), that if $\mu_a$ has a local \ho exponent equal to $\alpha$ at a point $t_0$, then $Z_a$ has at $t_0$ a pointwise \ho exponent equal to $\alpha^{1/H_a}-\log_3(4a-1)$. Hence the multifractal spectrum of $Z_a$ is deduced from the one of $\mu_a$ by the formula
$$d_{Z_a}(h) = \tilde d_{\mu_a} \left( (h+\log_3(4a-1))^{1/H}\right)$$
for every $h \in [(-\log_3 (a))^{1/H_a})-\log_3(4a-1), (-\log_3(2a-1))^{1/H_a})-\log_3(4a-1)]$. A more explicit formula is obtained as follows: for every $q\in\R$, if $\alpha = \tau_{\mu_a}$, then $\tilde d_{\mu_a} (\alpha) = q (\tau_{\mu_a})'(q) - \tau_{\mu_a}(q)$. 

The singularity spectra of $f_{2/3}$ and $\mu_{2/3}$ are given in Figure \ref{fig2}.

\smallskip

Finally, remark that the maximum of the spectrum is obtained for $\alpha_a=((\tau_{\mu_a})'(0))^{1/H_a}-\log_3(4a-1)$, and $d_{Z_a}(\alpha_a) =1$. After computations, we find $\alpha_a= \frac{-1}{3} (  \log_3 ( a^2(2a-1))$.

Let us consider the value of $a_0$ such that $\alpha_{a_0} =1$. Then $a_0^2(2a_0-1)=1/27$, i.e. $54a_0^3-27a_0^2 =1$.  When $a>a_0$, the set of points $t$ for which $h_{Z_a}(t)>1$ is of Lebesgue measure 1, hence we recover the main result of \cite{OKA}: $Z_a$ is differentiable on a set of Lebesgue measure 1. Here we obtain in addition the whole multifractal spectrum of $Z_a$. 



\subsection*{Acknowledgment}

The author thanks Yanick Heurteaux for a discussion on the oscillations of the self-similar functions.

\bibliographystyle{plain}

   \end{document}